\documentclass[a4paper,12pt]{amsart}
\usepackage[centering,text={15cm,23cm}]{geometry}
\usepackage{amscd}
\usepackage{amssymb}
\usepackage{csquotes}
\usepackage[foot]{amsaddr}
\usepackage{graphicx}
\usepackage{color}
\usepackage[all]{xy}
\usepackage{mathrsfs}
\usepackage{marvosym}
\usepackage{srcltx}
\usepackage{hyperref}

\definecolor{shadecolor}{rgb}{1,0.9,0.7}

\usepackage{amsfonts}
\usepackage{amsmath}
\usepackage{amscd}
\usepackage[T1]{fontenc}   
\usepackage[utf8]{inputenc} 
\usepackage{graphicx}       
\usepackage{amsmath}        
\usepackage{amssymb}        
\usepackage{amsthm}     
\usepackage{ae, aecompl}    
\usepackage{color}          
\usepackage{paralist}       
\usepackage{ifpdf}          
\usepackage{datetime}       
\usepackage{subfigure}      
\usepackage{wrapfig}        
\usepackage{lmodern}        
\usepackage{rotating}       
\usepackage{cancel}

\usepackage{seqsplit}
\usepackage{xstring}
\usepackage[final,color]{showkeys}        
\allowdisplaybreaks

\renewcommand{\showkeyslabelformat}[1]{
    \noexpandarg
    \StrSubstitute{#1}{ }{~}[\TEMP]
\fbox{\parbox{0.1\linewidth}{\normalfont\small\ttfamily\expandafter\seqsplit{\expandafter\TEMP}} 
}}        

\definecolor{refkey}{rgb}{0.6,0.7,0.6}
\definecolor{labelkey}{rgb}{0.6, 0.7, 0.6}

\usepackage{verbatim}         

\usepackage{mathtools}        

\usepackage{xifthen}

\usepackage{mathrsfs}       

   \usepackage{tikz}
   \usetikzlibrary{matrix}

\usepackage{rotating}    
\usepackage{listings}             

\newcommand{\SetDefine}[2]{
  \ifthenelse{\isempty{#2}}
             {\left\{ #1 \right\} }
{\left\{\vphantom{#2}#1\right.\left|\vphantom{#1}#2\right\}}
}

\newcommand{\tensor}{\otimes}
\newcommand{\loga}{\dagger}
\newcommand{\acknow}{ACKNOWLEDGEMENTS.}

\theoremstyle{plain} 
\newtheorem{theorem}{Theorem}[section]
\newtheorem{corollary}[theorem]{Corollary}
\newtheorem{lemma}[theorem]{Lemma}
\newtheorem{proposition}[theorem]{Proposition}

\theoremstyle{definition} 
\newtheorem{definition}[theorem]{Definition}
\newtheorem{example}[theorem]{Example}

\newtheorem{defn-prop}[theorem]{Definition-Proposition}
\newtheorem{remark}[theorem]{Remark}

\numberwithin{table}{section}
\numberwithin{figure}{section}

\DeclareMathOperator*{\fprods}{\sideset{}{_S}\prod}
\DeclareMathOperator{\Ima}{Im}
\DeclareMathOperator*{\fprodq}{\sideset{}{_{\mathfrak{Z}}}\prod}
\DeclareMathOperator{\Specan}{Spec_{an}}
\DeclareMathOperator{\Id}{Id}

\newcommand{\CC} {\mathbb{C}}

\newcommand {\foX}  {\mathfrak{X}}

\newcommand {\an}  {\mathrm{an}}

\newcommand {\coker} {\operatorname{coker}}

\newcommand {\Ex}  {\operatorname{Ex}}

\newcommand {\gp}  {{\operatorname{gp}}}

\newcommand {\Hom}  {\operatorname{Hom}}

\newcommand {\id}  {\operatorname{id}}

\renewcommand {\ker } {\operatorname{ker}}

\newcommand {\lra}  {\longrightarrow}

\newcommand {\M} {\mathcal{M}}

\renewcommand {\max} {{\operatorname{max}}}

\newcommand {\Mor} {\operatorname{Mor}}

\renewcommand{\O}  {\mathcal{O}}

\newcommand {\ol} {\overline}

\newcommand {\Spec} {\operatorname{Spec}}

\newcommand {\ul} {\underline}

\newcommand {\X} {\mathfrak X}

\begin{document}
\title[Analytic semi-universal deformations]{Analytic semi-universal deformations in logarithmic complex geometry}


\author{Raffaele Caputo}\address{Fachbereich Mathematik,
	Universit\"at Hamburg, Bundesstra\ss e~55, 20146 Hamburg,
	Germany}
\email{araffaelecaputo@gmail.com}



\maketitle

	\begin{abstract}
We show that every compact complex analytic space endowed with a fine logarithmic structure and every morphism between such spaces admit a semi-universal deformation. These results generalize the analogous results in complex analytic geometry first independently proved by A.~Douady and H.~Grauert in the '70. We follow Douady's two steps process approach consisting of an infinite-dimensional construction of the deformation space followed by a finite-dimensional reduction.    
	\end{abstract}

	\tableofcontents

	\section*{Introduction}
The main aim of this work is to extend the following  classical results in analytic deformation theory, Theorems \ref{thm:douady} and \ref{Flenner},  to the category of compact fine log complex spaces. 
\begin{theorem}\label{thm:douady}(Douady \cite{DouVers}, Grauert \cite{Gra}, Palamodov \cite{PalVers}, Forster--Knorr \cite{ForKno})
	Every compact complex analytic space admits a semi-universal deformation.
\end{theorem} 
\begin{theorem}(\cite[p. 130]{Flenner})\label{Flenner}
	Every morphism  between compact complex analytic spaces admits a semi-universal deformation.
\end{theorem}

We start by briefly reviewing some results in analytic deformation theory and by fixing some notation. For background material on complex analytic geometry, we recommend \cite{Fischer1976}, whereas, as references for log geometry, we recommend \cite{GrossKansas}, \cite{Ogus18} and \cite{KatNak}. The latter, in particular, explicitly deals with log structures on complex analytic spaces.
\begin{definition}
	Let $X_0$ be a compact complex analytic space. A deformation of $X_0$ is a triple $((S,0), X, i)$ consisting of a flat and proper morphism of complex spaces $\pi:X\to S$ and an isomorphism $i:X_0\to X(0)$, where $X(0):=\pi^{-1}(0)$. 
\end{definition}

 A deformation $((S,0), X, i)$  of a compact complex space $X_0$ is called \textit{complete}, if it contains, in a small neighborhood of the base point $0\in S$, all possible deformations of $X_0$.  Technically, this means that if $((T,0),Y,j)$ is another deformation of $X_0$, then there exists a morphism of germs $\varphi:(T,0)\to(S,0)$ and an isomorphism $\alpha:Y\to\varphi^{*}X$, such that $\alpha \circ j=\varphi^{*}i$. 
 
 Let $D:=(\{\cdot\},\mathbb{C}[\epsilon]/\epsilon^2)$ be the double point and $(S,0)$ a germ of complex spaces. Denote with $\Hom(D,(S,0))$ the set of morphisms of germs $D\to(S,0)$. We have a bijection
 \[\Hom(D,(S,0))\to \text{T}_0S\]
 sending $u:D\to (S,0)$ to $du(v)\in \text{T}_0S$,
 where $v\in \text{T}D$ is a basis element. If we denote with $\Ex^{1}(0)$ the set of isomorphism classes of  deformations of $X_0$ over $D$, we get a natural morphism
\begin{equation}\label{ks}
\text{ks}:\text{T}_0S\to \Ex^{1}(0),
\end{equation}
via $u\mapsto u^*\pi$. This morphism is called the \textit{Kodaira-Spencer map}. If $((S,0),X,i)$ is a complete deformation of $X_0$, then $\text{ks}$ is an epimorphism. If $\text{ks}$ is an isomorphism, the deformation is called \textit{effective} (see, for instance, the discussion in \cite[pp. 130--134]{Pal}).
In 1958, Kodaira, Nirenberg and Spencer (\cite{KodNirSpr}) proved that if $X_0$ is a compact complex manifold with $H^2(X_0;\mathcal{T}_{X_0})=~0$, then $X_0$ admits a complete and effective deformation with smooth base space. In 1962, Kuranishi (\cite{Kur}) proved the existence of a complete and effective deformation without the condition $H^2(X_0;\mathcal{T}_{X_0})=0$. In this case, the base space is a germ of complex spaces, in general singular. 
In 1964, A.~Douady (\cite{DouVar}), using his theory of Banach analytic spaces,  succeeded in giving a very elegant exposition of the results of Kuranishi. 
\begin{definition}(\cite[p. 601, Proposition 1]{DouVers}, \cite[p. 5, Definition 0.8]{Stie})\label{def:semiuniversal}
	Let $X_0$ be a compact complex space. A  deformation $((S,0),X,i)$ is called \textit{versal} if given any other deformation $((T,0), Y, j)$ of $X_0$, a subgerm $(T',0)$ of $(T,0)$ and a morphism $h':(T',0)\to (S,0)$ such that $Y|_{T'}\simeq h^{'*}X$,  there exists a morphism $h:(T,0)\to (S,0)$ such that $Y\simeq h^{*}X$ and $h|_{T'}=h'$.
\end{definition}
In literature, a versal and effective deformation is called \textit{semi-universal} or \textit{miniversal}. By a general result  of H.~Flenner (\cite[Satz 5.2]{Flenner1981}),  every versal deformation gives a semi-universal deformation.
\smallskip

We outline the key ideas in Douady's construction of a semi-universal deformation of a compact complex space. We start by noticing that we can cover  a compact complex space $X_0$  with finitely many open subsets $(U_i)_{i\in I_0}$, such that, for each $i\in I_0$,
there exists a closed subset $Z_i\subset W_i$, for some $W_i$ open in $\mathbb{C}^{n_i}$, and an isomorphism
\begin{equation}\label{morph_dis}
f_i:Z_i\to U_i.
\end{equation}
Moreover, we can find an isomorphism of the form \eqref{morph_dis} for any double $U_{ij}:=U_i\cap U_j$ and triple $U_{ijk}:=U_i\cap U_j\cap U_k$ intersection.
 The collection of closed subspaces $((Z_i),(Z_{ij}),(Z_{ijk}))$ is a \textit{disassembly} of $X_0$, where the \textit{assembly instructions} are encoded into the isomorphisms $((f_i),(f_{ij}),(f_{ijk}))$ via the transition maps $(f_i^{-1}\circ f_j)$.
A deformation of $X_0$ is obtained by deforming each closed subspace $Z_i$, together with the gluing morphisms $f_i$, and by assembling together the obtained deformed subspaces.

Douady's key insight was to choose  special (\enquote{\textit{privileged}}) subspaces $(Y_i)$ of given polycylinders $(K_i\subset\mathbb{C}^{n_i})$ for the closed subspaces $(Z_i)$, and to show that the collection of all privileged  subspaces of a given polycylinder can be endowed with an analytic structure.
More precisely, given a polycylinder $K_i\subset\mathbb{C}^{n_i}$, we can consider the Banach algebra
\[B(K):=\{h:K_i\to\mathbb{C}| h\text{ is continuous on }K_i\text{ and analytic on its interior}\}.\]   
An ideal $I\subset B(K)$ is called \textit{direct} if there exists a $\mathbb{C}$-vector subspace $J$ of $B(K)$, such that $B(K)=I\oplus J$ as $\mathbb{C}$-vector spaces. Douady showed in \cite[p. 34]{DouThesis}, that the set
\[\mathcal{G}(B(K)):=\{I\subset B(K)| I \text{ is direct}\}\]
can be endowed with the structure of a Banach manifold (see \cite[p. 16]{DouThesis}; \cite[p. 38, Example 3.15]{Raf}). The space $\mathcal{G}(B(K))$ is called the \textit{Grassmannian} of $B(K)$. Furthermore, if we consider $B(K)$ as a module over itself, the set
\[\mathcal{G}_{B(K)}(B(K)):=\{I\in \mathcal{G}(B(K))| I \text{ is a }B(K)\text{-submodule of } B(K)\}\]
can be endowed with the structure of a Banach analytic space (see \cite[pp. 29--30]{DouThesis}; \cite[p. 39, Example 3.21]{Raf}) and the subset 
 \begin{equation}\label{grassmannian}
 \mathcal{G}(K):=\{I\in\mathcal{G}_{B(K)}(B(K))|I \text{ admits a finite free resolution}\}
 \end{equation}   is  open in $\mathcal{G}_{B(K)}(B(K))$. The \textit{privileged} subspaces of a given polycylinder $K_i$ are precisely those subspaces corresponding to the \textit{direct ideals of $B(K)$ admitting a finite free resolution} (see \cite[p. 577]{DouVers} and \cite[p. 256]{Lepot}). 
 In \cite[p. 62, Theorem 1]{DouThesis}, Douady showed that every compact complex space can be covered with finitely many privileged subspaces of polycylinders.
 
Now, given a covering of a compact complex space $X_0$ with \textit{privileged} charts $(f_i:Y_i\to X_0)$, since intersections of privileged polycylinders are not in general privileged, one needs to cover the intersections too. In order to have the transition maps well-defined, one needs to work with two polycylinders 
\begin{equation}\label{doublesubset2}
\tilde{K}_i\subset \mathring{K}_i
\end{equation}
 for double intersections and three polycylinders 
\begin{equation}\label{doublesubset}
K'_i\subset \mathring{\tilde{K}}_i, \tilde{K}_i\subset \mathring{K}_i
\end{equation} 
for triple intersections. We rewrite \eqref{doublesubset2} and \eqref{doublesubset} using the following notation
\begin{equation}\label{doublesubsetnotation2}
\tilde{K}_i\Subset K_i
\end{equation} 
and
\begin{equation}\label{doublesubsetnotation}
K'_i\Subset \tilde{K}_i\Subset K_i
\end{equation} 
respectively. Let
\begin{equation}\label{type cuirasse}
\mathfrak{I}:=(I_{\bullet}, (K_{i})_{i\in I},(\tilde{K}_{i})_{i\in I},(K'_{i})_{i\in I_{0}\cup I_{1}}),
\end{equation}
where  $I_{\bullet}$ is a finite simplicial set of dimension 2 (see, for instance, \cite[p. 587]{DouVers}) and the collections of polycylinders satisfy \eqref{doublesubsetnotation2} and \eqref{doublesubsetnotation}.
A \textit{cuirasse} $q$ of \textit{type} $\mathfrak{I}$ on a compact complex space $X_0$ is a \textit{disassembly} of $X_0$ given by a collection of pairs $q:=\{(Y_i,f_i)\}_{i\in I}$, where $Y_i\subset K_i$ is privileged, $f_i:Y_i\to X_0$ is a morphism, and they satisfy gluing relations on double and triple intersections (see \cite[p. 587]{DouVers}).

In \cite[p. 588]{DouVers}, Douady showed that the set of all cuirasses of a fixed type $\mathfrak{I}$ on a compact complex space $X_0$
\begin{equation}\label{cuirasses}
\mathcal{Q}(X_0):=\{q \text{ is a cuirasse on } X_0\}
\end{equation}
can be endowed with the structure of a Banach analytic space. Moreover, if $X\to S$ is a deformation of $X_0$, a choice of a cuirasse $q_s$ on each fibre $X_s$ is called a \textit{relative cuirasse} on $X$ over $S$. More precisely, in \cite[p. 588]{DouVers}, Douady showed that the set
\begin{equation}\label{space of rel cuirasses}
\mathcal{Q}_{S}(X):=\{(s,(Y_{i},f_{i})_{i\in I})\lvert s\in S, (Y_{i},f_{i})_{i\in I}\in \mathcal{Q}(X(s))\},
\end{equation}
that is
\[\mathcal{Q}_{S}(X)=\bigsqcup_{s\in S}\mathcal{Q}(X(s)),\]
can be endowed with the structure of a  Banach analytic space. Then, a \textit{(local) relative cuirasse}  on $X$ over $S$ is defined as a (local) section
\begin{equation}\label{rel cuirasses}
q:S\to \mathcal{Q}_{S}(X)
\end{equation}
of the natural projection $\pi: \mathcal{Q}_{S}(X)\to S$.

On the other side stands the notion of \textit{puzzle}. Informally speaking, a  puzzle is a compact complex space delivered in pieces, together with the assembly manual. Technically, a puzzle $z$ is given by a collection $z:=\{(Y_i,g^j_i)\}_{i\in I, j\in\partial i}$, where $Y_i\subset K_i$ is a privileged subspace and $g^j_i:Y_j\to Y_i$ is a morphism. This collection of data satisfies gluing axioms (\cite[p. 589]{DouVers}). The collection of puzzles 
\begin{equation}\label{puzzles}
\mathfrak{Z}:=\{(Y_i,g^j_i)_{i\in I, j\in\partial i}\}
\end{equation}
form a Banach analytic space, each puzzle $z$ glues to a compact complex space $\mathfrak{X}_z$  and the collection of compact complex spaces $(\mathfrak{X}_z)_{z\in\mathfrak{Z}}$ glues to a proper Banach analytic family $\mathfrak{X}$ over $\mathfrak{Z}$ (see \cite[p. 591]{DouVers}), which is anaflat (see \cite[p. 66, Definition and Proposition 1]{DouThesis}). 

Now, let $X\to S$ be a deformation of $X_0$. The aim is to produce a map $\varphi:S\to \mathfrak{Z}$, such that, in a neighborhood of some base point $z_0\in\mathfrak{Z}$, with $\mathfrak{X}_{z_0}\simeq X_0$, we have $\varphi^{*}\mathfrak{X}\simeq X$. To achieve this end, a special role is played by \textit{triangularly privileged cuirasses} on $X_0$ (see \cite[p. 588]{DouVers}). Informally speaking, these are cuirasses on $X_0$ that extend to cuirasses on the nearby fibres $X_s$. Douady showed that every compact complex space $X_0$ admits a \textit{triangularly privileged cuirasse} (\cite[p. 588]{DouVers}) 
\begin{equation}\label{triang}
q_0\in \mathcal{Q}(X_0).
\end{equation}
This means that if $X\to S$ is a deformation of $X_0$, with base point $0\in S$, and $q_0$ is a triangularly privileged cuirasse on $X_0$, then we get the existence of  a continuous family of cuirasses $\{q_s\}_{s\in S}$, where $q_s$ is a cuirasse on the fibre $X_s$, for $s$ in a small neighborhood of $0$. Namely, we can find a  \textit{(local) relative cuirasse} $q:S\to\mathcal{Q}_S(X)$ on $X$ over $S$, such that $q(0)=q_0$. 
Now, since every cuirasse $q_s=\{(Y_i,f_i)\}$ naturally produces an \textit{associated} puzzle (\cite[p. 590]{DouVers}) via
\begin{equation}\label{prop:puzzle associated}
z_{q_s}:=(Y_i,g_{i}^{j}:=f_{i}^{-1}\circ f_{j})_{i\in I,j\in\partial i},
\end{equation}
we get a morphism (\cite[p. 591]{DouVers})
 \begin{equation}\label{def:morph associated to a cuirasse}
 \begin{split}
 \varphi_q:S&\to\mathfrak{Z}\\
s&\mapsto z_{q_s}.\\
 \end{split}
\end{equation}

Because a cuirasse $q_s$ is a disassembly of a compact complex space $X_s$ and the associated  puzzle $z_{q_s}$ glues to a compact complex space $\mathfrak{X}_{z_{q_s}}$, it is reasonable to expect that $\mathfrak{X}_{z_{q_s}}$ is isomorphic to $X_s$. In fact, we have an $S$-isomorphism (\cite[p. 592]{DouVers})
\begin{equation}\label{prop:isomorphism alpha}
\alpha_q:\varphi^{*}_q\mathfrak{X}\to X.
\end{equation}
In other words, the Banach analytic family $\mathfrak{X}\to \mathfrak{Z}$ contains all possible deformations of $X_0$ in a neighborhood of $z_{q_0}$. That is, the family is \textit{complete}.

 An involved finite-dimensional reduction procedure (\enquote{\textit{a cure d'amaigrissement}}) is  used to obtain a finite-dimensional semi-universal deformation of $X_0$ out of the complete infinite-dimensional family $\mathfrak{X}\to \mathfrak{Z}$ (see \cite[pp. 593--599]{DouVers}, \cite[pp. 20--46]{Stie} and subsection \ref{finite dimensional reduction}).  This ends our survey about Douady's construction of  a semi-universal deformation of a compact complex space. 
\medskip

Now, we assume that $X_0$ comes endowed with a fine log structure $\mathcal{M}_{X_0}$. We view $X_0$ as a log space over the point $\Spec\CC$ with trivial log structure.

\begin{definition}\label{def:versal def log spaces}
	A \textit{deformation of a compact fine log complex space} $(X_0,\mathcal{M}_{X_0})$ is a triple $((S,s_0), (\mathfrak{X},\mathcal{M}_\mathfrak{X}), i)$, where $S$ is a complex space endowed with trivial log structure, $s_0\in S$,
	 $p:(\mathfrak{X},\mathcal{M}_\mathfrak{X})\rightarrow (S,\mathcal{O}_S^{\times})$  is a log morphism between fine log complex spaces with underlying map of complex spaces $\foX\to S$ proper and flat, and $i:(X_0,\mathcal{M}_{X_0})\rightarrow (\mathfrak{X},\mathcal{M}_\mathfrak{X})(s_0):=p^{-1}(s_0)$ is a log isomorphism.
\end{definition}
A deformation is \textit{complete} if for any other deformation $((T,t_0),(X,\mathcal{M}_X),j)$ of $(X_0,\mathcal{M}_{X_0})$, there exists a morphism
$\psi:(T,\mathcal{O}_T^{\times})\rightarrow (S,\mathcal{O}_S^{\times})$, sending $t_0$ to $s_0$, and a log $T$-isomorphism \[\alpha:(X,\mathcal{M}_X)\rightarrow (\mathfrak{X},\mathcal{M}_\mathfrak{X})\times_{(S,\mathcal{O}_S^{\times})}(T,\mathcal{O}_T^{\times}),\] such that $\alpha\circ j=i$.
For the sake of readability, in what follows, we shall mostly denote a complex space endowed with trivial log structure $(S,\O^{\times}_S)$ just by $S$.

One of the key points, in the construction of deformations of log spaces, is to find a proper way to deform the log structure $\mathcal{M}_{X_0}$ \textit{coherently} with the deformation of the underlying analytic space $X_0$. We show, in subsection \ref{sub:gluing charts}, that we can \textit{disassemble} $\mathcal{M}_{X_0}$ using \textit{log charts} satisfying gluing conditions on double and triple intersections (Proposition \ref{prop:existence directed log charts}).  That is, the log structures associated to the log charts glue to a global log structure $\mathcal{M}^a_{X_0}$ on $X_0$ isomorphic to $\mathcal{M}_{X_0}$. We call this collection of log charts a set of \textit{directed} log charts (Definition \ref{def:compatible log charts}). This insight leads to the notion of \textit{log cuirasse} (Definition \ref{def:log cuirasse}) and \textit{log puzzle} (Definition \ref{def:log puzzle}).

In subsection \ref{Infinite dimensional construction}, we construct an infinite-dimensional \textit{log} family $(\mathfrak{X},\mathcal{M}_{\mathfrak{X}})\to\mathfrak{Z}^{\log}$ (Proposition \ref{prop:log moduli}).
Given  a log deformation $(\mathcal{Y},\mathcal{M}_{\mathcal{Y}})\to T$ of $(X_0,\mathcal{M}_{X_0})$, with base point $t_0$,
an essential point is to show that a triangularly privileged \textit{log} cuirasse $q_0^{\dagger}$ exists on $(X_0,\mathcal{M}_{X_0})\simeq (\mathcal{Y},\mathcal{M}_{\mathcal{Y}})(t_0)$ and it extends to a log cuirasse $q^{\dagger}_t$ on the fibre $(\mathcal{Y}_t,\mathcal{M}_{\mathcal{Y}_t})$, for $t$ in a neighborhood of $t_0$ (Propositions \ref{prop:log fibrewise cuirasse} and \ref{prop:log triangularly smoothness}). This allows us to show the \textit{completeness} of the log family $(\mathfrak{X},\mathcal{M}_{\mathfrak{X}})\to\mathfrak{Z}^{\log}$. 

In subsection \ref{finite dimensional reduction},
we proceed with a finite-dimensional reduction procedure, which produces a semi-universal deformation of $(X_0,\mathcal{M}_{X_0})$ out of the complete log family $(\mathfrak{X},\mathcal{M}_{\mathfrak{X}})\to\mathfrak{Z}^{\log}$. The finite-dimensionality is achieved with the exact same procedure used by Douady in the classical case.  This is because the space $\mathfrak{Z}^{\log}$ of log puzzles does not come endowed with a non-trivial log structure. 
We prove
\begin{theorem}\label{thm 1}(Theorem $\ref{thm:log fin dim douady versal}$)
	Every  compact fine log complex space $(X_0,\mathcal{M}_{X_0})$ admits a semi-universal deformation $((S,s_0), (\mathfrak{X},\mathcal{M}_\mathfrak{X}), i)$.
\end{theorem}
For a construction of a semi-universal deformation in the \textit{non-fine} log context see, for instance, \cite{RS} where a semi-universal family is obtained by means of \textit{Artin approximation} (see, also, \cite{ruddat2019local}).
\medskip

 The existence of semi-universal deformations of morphisms between compact complex analytic spaces follows naturally from Douady's results (see \cite[p. 130]{Flenner}). Analogously, we take a further step in our work studying semi-universal deformations of log morphisms. 
 Given a morphism of log complex spaces, we have the notion of \textit{log smoothness} (see, for instance, \cite[p. 107]{GrossKansas}) and \textit{log flatness} (see \cite{illusie2013}). These notions  generalize and extend the classical notions of smoothness and flatness, which are retrieved if we consider complex spaces endowed with trivial log structures. In \cite{KatoMa}, K. Kato writes
that a log structure is \enquote{magic by which a degenerate scheme begins to behave as being
non-degenerate}.

For example, the affine toric variety $\Spec_{\an}\mathbb{C}[P]$,
with its canonical divisorial log structure, is log smooth over $\Spec\mathbb{C}$ (equipped with the trivial log
structure), despite almost always not being smooth in the usual sense.
 In what follows, we denote the analytic spectrum  $\text{Spec}_{\text{an}}\mathbb{C}[P]$ of a monoid ring simply by $\Spec\mathbb{C}[P]$.

 \smallskip
 
In section \ref{log morph}, we prove the following

\begin{theorem}\label{thm 2}(Theorem $\ref{thm:vers log morph}$ and Proposition $\ref{prop:log flatness}$)
Every morphism of compact fine log complex spaces $f_0:(X_0,\mathcal{M}_{X_0})\rightarrow (Y_0,\mathcal{M}_{Y_0})$ admits a semi-universal deformation $f$ over a germ of complex spaces $(S,s_0)$. Moreover, if $f_0$ is log flat (or log smooth), then  $f$ is log flat (or log smooth) in an open neighborhood of $s_0$.
\end{theorem}

As a corallary result (Corollary \ref{cor:def log morph}), we obtain a \textit{relative} semi--universal deformation of a compact fine log complex space  $(X_0,\mathcal{M}_{X_0})$ \textit{over} a fine log complex space $(Y_0,\mathcal{M}_{Y_0})$ (Definition \ref{def:relative 1}). Notice that, in this case,   $Y_0$ needs not to be compact. If $(X_0,\mathcal{M}_{X_0})$ is a log subspace of $(Y_0,\mathcal{M}_{Y_0})$, we get a semi--universal deformation of a log \textit{subspace} in a \textit{fixed ambient} log space (Remark \ref{rem:relative 3}).
\medskip

The focus of this work is the construction of analytic deformations via Douady's patching method rather than a comprehensive treatment of deformations of analytic log spaces. In particular, we do not discuss infinitesimal or formal deformations. The classical treatment of these topics in the algebraic geometric setup (see \cite{KatoFI} and \cite{KatoF})
readily carry over to the analytic setup treated here. See also \cite{Felten}, for a more recent treatment of log smooth deformations from the point of view of differential graded algebras.
\bigskip

\acknow{This work is part of my Ph.D. thesis written at the University of Hamburg under the supervision of Bernd Siebert, whom I thank for the help and inspiration which have made the realization of this work possible. 
I also thank Siegmund Kosarew for very helpful comments, Helge Ruddat and Simon Felten for the hospitality at the Johannes Gutenberg University of Mainz during spring 2019 and Mark Gross for the hospitality at the University of Cambridge during spring 2020. 
I thank Bernd Siebert and The University of Texas at Austin for financial support.
At the University of Hamburg, I was supported by the Research Training Group 1670 \enquote{Mathematics inspired by String Theory and Quantum Field Theory}, funded by the German Research Foundation -- Deutsche Forschungsgemeinschaft (DFG). }

	\section[Semi-universal def. of compact fine log complex spaces]{Semi-universal deformations of compact fine log complex spaces}

In  what follows, we construct a  semi-universal deformation in the general case of a compact complex space $X_{0}$ endowed with a fine log structure $\mathcal{M}_{X_{0}}$. 
	
\subsection{Gluing log charts}
\label{sub:gluing charts}
Let $(X_0,\M_{X_0})$ be a compact fine log complex space. Denote by $\alpha:\M_{X_0}\to\O_{X_0}$ the structure map. The sheaf of monoids
\[\overline{\mathcal{M}}_{X_0}:=\mathcal{M}_{X_0}/\Ima \alpha^{-1}=\mathcal{M}_{X_0}/\mathcal{O}^{\times}_{X_0},\]
written additively, is called the \textit{ghost sheaf} of  $\mathcal{M}_{X_0}$.	
We assume that $\overline{\mathcal{M}}^{\text{gp}}_{X_0}$ is torsion free.

\smallskip

We want to find a universal setup for constructing log structures from gluing of log
charts. This is quite analogous to the case of sheaves, see for example \cite[Exercise II.1.22]{Hart}. 
Assume we have a covering of $X_0$ by open sets $U_i$ for an ordered index set $J_0$, and for each $U_i$ a log chart 
\[
\theta_i:P_i\lra \Gamma(U_i,\M_{X_0}).
\]
We identify $\theta_i$ with the corresponding map of monoid sheaves $\ul P_i\to
\M_{X_0}|_{U_i}$. 
For $l=1,2$, set \[J_l:=\{(i_0,...,i_l)\in J_0^{l+1}: U_{i_0}\cap...\cap U_{i_l}\neq \emptyset\}.\]  
We get maps $d_m:J_l\rightarrow J_{l-1}$, for $0\leq m\leq l$ and $1\leq l\leq 2$, sending $(i_0,..,i_m,..,i_l)$ to $(i_0,..,i_{m-1},i_{m+1},..,i_l)$. We set 
\begin{alignat*}{9}
J &:=\bigcup^2_{l=0}J_l,\\
\partial{i} &:=\left\lbrace d_{0}i, ...., d_{l}i \right\rbrace && \text{, if } i\in J_{l}.
\end{alignat*}
The set $J$, together with the maps $(d_m)$, is called a \textit{simplicial set of order 2}. 

For each $j:=(i_0,i_1)\in J_1$, assume that there is a log chart
\[
\theta_{j}:P_{j}\lra \Gamma(U_{j},\M_{X_0})
\]
and comparison maps
\begin{eqnarray*}
	\varphi^{i}_{j}: P_{i}&\longrightarrow & P_{j}\oplus \Gamma(U_{j},\O^{\times}_{X_0}),
\end{eqnarray*}
for $i\in\partial j$, with the property
\begin{equation}
\label{Eqn: Comparison equation}
\big(\theta_{j}\cdot\Id_{\O^{\times}_{X_0}|_{U_j}}\big)\circ\varphi^{i}_{j}= \theta_i|_{U_{j}}.
\end{equation}

Each $\theta_i$ defines an isomorphism
of $\M_{U_i}$ with the log structure $\M_i$ associated to the pre-log structure $\beta_i:=\alpha\circ\theta_i$.
Similarly, the pre-log structure $\beta_{j}:=\alpha\circ\theta_{j}$ defines a log structure $\M_{j}$ and $\theta_j$ defines an
isomorphism of log structures $\M_{U_{j}}\simeq \M_{j}$. From this point of
view, equation~\eqref{Eqn: Comparison equation} means that $\varphi^i_{j}$
provides an isomorphism between $\M_i|_{U_{j}}$ and $\M_{j}$, and this
isomorphism is compatible with the isomorphisms $\M_i|_{U_{j}}\simeq
\M_{U_{j}}$ and $\M_{j}\simeq \M_{U_{j}}$. 

Now, if we have $\theta_i$, $\theta_{j}$, $\varphi^i_{j}$, fulfilling \eqref{Eqn: Comparison equation}, we need compatibility on triple intersections for the patching of the $\M_i$ to be consistent. To formulate this cocycle condition in terms of log charts, assume, for each $k:=(i_0,i_1,i_2)\in J_2$, a third system of charts
\[
\theta_{k}: P_{k}\lra\Gamma(U_{k},\M_{X_0})
\]
and comparison maps
\begin{eqnarray*}
	\varphi^{j}_{k}: P_{j}&\lra& P_{k}\oplus \Gamma(U_{k},\O^{\times}_{X_0}),
\end{eqnarray*}
for $j\in \partial k$. The analogue of the compatibility condition \eqref{Eqn: Comparison equation} is
\begin{equation}
\label{Eqn: triple compatibility}
\big(\theta_{k}\cdot\Id_{\O^{\times}_{X_0}|_{U_{k}}}\big)\circ\varphi^{j}_{k}= \theta_{j}|_{U_{k}}.
\end{equation}
Again, the $\varphi^{j}_{k}$ define an isomorphism between the log structure $\M_{j}|_{U_{k}}$ on $U_{k}$ and the log structure $\M_{k}$ associated to the pre-log structure $\beta_{k}:=\alpha\circ\theta_k$. In particular, all the isomorphisms of log structures 
are compatible and the $(\M_i)_{i\in J_0}$ glue in a well-defined fashion, as do their structure maps, to a log structure on $X_0$ isomorphic to $\M_{X_0}$. This is just standard sheaf theory, for sheaves of monoids.

\begin{definition}\label{def:compatible log charts}
	A set of \textit{directed} log charts is  a set of log charts $(\theta_i:P_i\to\mathcal{M}_{U_i})_{i\in J}$ covering $(X_{0},\mathcal{M}_{X_0})$, together with a morphism
	\[\varphi^{i}_j:P_i\to P_j\oplus \mathcal{O}^{\times}_{U_j},\]
	for each $j\in J_1\cup J_2$ and $i\in\partial j$, such that
	\[(\theta_j\cdot \Id_{\mathcal{O}^{\times}_{U_j}})\circ \varphi^i_j=\theta_i|_{U_j}.\]
\end{definition} 
In Proposition \ref{prop:existence directed log charts}, we show that every compact fine log complex space can be covered with a finite set of directed log charts.
\medskip
	
	Now, let us forget that the $(\theta_i)_{i\in J_0}$, $(\theta_{j})_{j\in J_1}$ and $(\theta_{k})_{k\in J_2}$ are charts for the given log structure. Let  $(U_i)_{i\in J_0}$ be an open cover of $X_0$ and $J$ as above. 
	Assume we have pre-log structures $(\beta_i)_{i\in J}$ and comparison maps $(\varphi^i_{j})_{j\in J_1\cup J_2,i\in\partial j}$ satisfying equations~\eqref{Eqn: Comparison equation} and \eqref{Eqn: triple compatibility}. Then the log structures $(\M_i)_{i\in J_0}$ glue to a log structure $\M_{X_0}$ on $X_0$ in such a way that the gluing data $(\beta_{j})_{j\in J_1}$ and compatibility $(\beta_{k})_{k\in J_2}$ arise from identifying $\M_{j}$ and $\M_{k}$ with restrictions of $\M_{X_0}$ to $U_{j}$ and $U_{k}$ respectively. 
	

	\begin{definition}\label{def:pre-log atlas}
		Let $X_0$ be a compact complex space. With the above notation, we call a \textit{pre-log atlas} on $X_0$ a collection of data \[\{(\beta_i:P_i\rightarrow\mathcal{O}_{U_i})_{i\in J},(	\varphi^{i}_{j}: P_{i}\to P_{j}\oplus \O^{\times}_{U_j})_{j\in J_1\cup J_2,i\in\partial j}\}\] satisfying
		\begin{equation}
		\big(\beta_{j}\cdot\Id_{\O^{\times}_{U_j}}\big)\circ\varphi^{i}_{j}= \beta_i|_{U_{j}}.
		\end{equation}
	\end{definition}	
	\subsection{Infinite dimensional construction}
	\label{Infinite dimensional construction}
	
The notion of \text{log structure} can be naturally extended to the category of Banach analytic spaces. Indeed, let $(X,\Phi)$ be a Banach analytic space (see \cite[pp. 22--25]{DouThesis}; \cite[p. 38, Definition 3.16]{Raf}). Setting $\O_X:=\Phi(\mathbb{C})$, we get a ringed space $(X,\mathcal{O}_X)$. 
\begin{definition}
	A \textit{pre-log structure} on a Banach analytic space $(X,\Phi)$ is a sheaf of monoids $\mathcal{M}_{X}$ on $X$ together with a homomorphism of sheaves of monoids:
	\[\alpha_{X}:\mathcal{M}_{X}\rightarrow\mathcal{O}_{X},\]
	where the monoid structure on $\mathcal{O}_X$ is given by multiplication.
	A pre-log structure is a called a \textit{log structure} if
	\[\alpha_{X}:\alpha_{X}^{-1}(\mathcal{O}_{X}^{\times})\rightarrow\mathcal{O}_{X}^{\times}\]
	is an isomorphism. 
\end{definition}
The notion of \textit{fine} log structure extends naturally to the Banach analytic setting. 
In what follows, we shall mostly denote 
a log Banach analytic space endowed with the trivial log structure $(S,\O^{\times}_S)$ just by $S$. Moreover, for the sake of readability, we shall often write
Banach analytic morphisms just set-theoretically.
\smallskip

Let  $(X_0, \mathcal{M}_{X_{0}})$  be a compact fine log complex space.
	
	\begin{definition}\label{def:discrete part}
		 Let $(\theta_{i,0}:P_i\rightarrow \mathcal{M}_{U_i})_{i\in J}$, with comparison morphisms \[(\varphi^i_{j,0}:=(\phi^{i}_{j,0},\eta^i_{j,0}):P_i\to P_j\oplus\mathcal{O}^{\times}_{U_j})_{j\in J_1\cup J_2, i\in\partial j},\]  be a finite set of directed log charts covering $(X_0,\M_{X_0})$ (Definition $\ref{def:compatible log charts}$) given by Proposition $\ref{prop:existence directed log charts}$.
	\end{definition}
	Let
	\[\mathfrak{I}:=(I_{\bullet}, (K_{i})_{i\in I},(\tilde{K}_{i})_{i\in I},(K'_{i})_{i\in I_{0}\cup I_{1}})\]
	be as in $\eqref{type cuirasse}$ and $\mathfrak{Z}$ the space of puzzles $\eqref{puzzles}$.	
	Without loss of generality, we assume that the index sets $I$ and $J$ (Definition $\ref{def:discrete part}$) coincide. 	 We can define the notion of \textit{log puzzle}, which, informally speaking, is a compact fine log complex space delivered in pieces with the instructions to glue them together.
	\begin{definition}\label{def:log puzzle}
		
		A \textit{log puzzle} is a pair $(z,l)$, where $z:=(Y_{i},g^{i}_j)\in\mathfrak{Z}$ is a puzzle and $l$ is a collection of data 
		\[((\beta_i:P_i\rightarrow\mathcal{O}_{Y^{\circ}_i})_{i\in I},(\eta^i_{j}:P_{i}\to  \O^{\times}_{Y^{\circ}_j})_{j\in I_1\cup I_2,i\in\partial j}),\]
		satisfying
		\begin{equation}
		\big((\beta_{j}\cdot\Id_{\O^{\times}_{Y^{\circ}_{j}}})\circ\varphi^{i}_{j}= \beta_i|_{Y^{\circ}_{j}}\big)_{j\in I_1\cup I_2, i\in\partial j},
		\end{equation}
		where $\varphi^{i}_{j}:=(\phi^{i}_{j,0},\eta^i_{j})$, with $\phi^{i}_{j,0}:P_i\to P_j$ given by Definition $\ref{def:discrete part}$.
	\end{definition}

	\begin{definition}\label{def:log puzzle set}
		We denote the set of log puzzles by $\mathfrak{Z}^{\log}$.
	\end{definition}
	
	
	The set of log puzzles $\mathfrak{Z}^{\log}$ can be endowed with a Banach analytic structure. Indeed, for each polycylinder $K_i$, let us consider the Grassmannian $\mathcal{G}(K_i)$ (\ref{grassmannian}) and let $\Id:\mathcal{G}(K_i)\to\mathcal{G}(K_i)$ be the identity map. Identifying $\Id$ with its graph, we get a universal $\mathcal{G}(K_i)$-anaflat subspace $\underline{Y}_i\subset \mathcal{G}(K_i)\times\mathcal{G}(K_i)\subset\mathcal{G}(K_i)\times K_i$ (see \cite[p. 579]{DouVers}, \cite[pp. 258--259]{Lepot} and \cite[p. 183, Theorem 4.13]{Pourpriv}). 
	Let us consider the Banach analytic space
	\[\mathfrak{M}:=\fprodq_{i\in I} \mathfrak{Mor}_{\mathfrak{Z}}(\underline{Y}_i,\Spec\mathbb{C}[P_{i}]\times \mathfrak{Z})\times_{\mathfrak{Z}}\fprodq_{j\in I_1\cup I_2}\fprodq_{i\in\partial j}\mathfrak{Mor}_{\mathfrak{Z}}(\underline{Y}_j,\Spec\mathbb{C}[P^{\gp}_{i}]\times \mathfrak{Z}).\]
	Each morphism \[\eta^{i}_j(z)^{\gp}:P^{\gp}_i\to\mathcal{O}^{\times}_{\underline{Y}^{\circ}_{j}(z)}\] induces a morphism $\eta^i_j(z):P_i\to\mathcal{O}^{\times}_{\underline{Y}^{\circ}_{j}(z)}$. Hence, a point in $\mathfrak{M}$ can be written as
	\[(z,(\beta_{i}(z):P_i\to\mathcal{O}_{\underline{Y}^{\circ}_{i}(z)})_{i\in I},(\eta^i_j(z):P_i\to\mathcal{O}^{\times}_{\underline{Y}^{\circ}_{j}(z)})_{j\in I_1\cup I_2,i\in\partial j}).\]
	Thus, we naturally get an injective map
	
	\begin{equation}\label{eq:injective map}
	\rho:\mathfrak{Z}^{\log}\hookrightarrow \fprodq_{i\in I} \mathfrak{Mor}_{\mathfrak{Z}}(\underline{Y}_i,\Spec\mathbb{C}[P_{i}]\times \mathfrak{Z})\times_{\mathfrak{Z}}\fprodq_{j\in I_1\cup I_2}\fprodq_{i\in\partial j}\mathfrak{Mor}_{\mathfrak{Z}}(\underline{Y}_j,\Spec\mathbb{C}[P^{\text{gp}}_{i}]\times \mathfrak{Z}).
	\end{equation}

	\begin{proposition}
		The universal space of log puzzles $\mathfrak{Z}^{\log}$ is Banach analytic.
	\end{proposition}
	\proof
	Let $(\phi_{j,0}^i:P_i\to P_j)_{j\in I_1\cup I_2,i\in \partial j}$ given by Definition $\ref{def:discrete part}$, we set \[\varphi^{i}_{j}(z):=(\phi_{j,0}^i,\eta^i_j(z)):P_i\to P_j\oplus \mathcal{O}^{\times}_{\underline{Y}^{\circ}_{j}(z)}.\]
	The subset $\rho(\mathfrak{Z}^{\log})$ of $\mathfrak{M}$ ($\ref{eq:injective map}$) is
	defined by the equations
	\begin{equation}
	\big((\beta_{j}(z)\cdot\Id_{\O^{\times}_{\underline{Y}^{\circ}_{j}(z)}})\circ\varphi^{i}_{j}(z)= \beta_i(z)|_{\underline{Y}^{\circ}_{j}(z)}\big)_{j\in I_1\cup I_2, i\in\partial j}.
	\end{equation}
	Thus, we can define a double arrow
	\[(\rho_1,\rho_2):\mathfrak{M}
	\rightrightarrows\fprodq_{j\in I_1\cup I_2}\mathfrak{Mor}_{\mathfrak{Z}}(\underline{Y}_{j},\Spec\mathbb{C}[P_{j}]\times \mathfrak{Z})\]
	by
	\[\rho_1:(z,(\beta_i(z)),(\eta^{i}_j(z)))\mapsto(z,(\beta_i(z)|_{\underline{Y}^{\circ}_{j}(z)}))\]
	and 
	\[\rho_2:(z,(\beta_i(z)),(\eta^{i}_j(z)))\mapsto(z,((\beta_{j}(z)\cdot\Id_{\O^{\times}_{\underline{Y}^{\circ}_{j}(z)}})\circ\varphi^{i}_{j}(z))).\]
	Then $\mathfrak{Z}^{\log}$ is given by the kernel of the double arrow defined by $\rho_1$ and $\rho_2$:
	\[\mathfrak{Z}^{\log}=\ker(\rho_1,\rho_2).\]

	\endproof	
	Let $p:\mathfrak{Z}^{\log} \rightarrow \mathfrak{Z}$ be the canonical projection and consider the Banach analytic space $\mathfrak{X}_{\log}:=p^*\mathfrak{X}$ over $\mathfrak{Z}^{\log}$.
	 
	\begin{proposition}\label{prop:log moduli}
		The Banach analytic space $\mathfrak{X}_{\log}$ comes naturally endowed with a fine log structure $\mathcal{M}_{\mathfrak{X}_{\log}}$.
	\end{proposition}
	\proof
	%
	By \cite[p. 192, Theorem 5.13]{Pourpriv} (see, also, \cite[p. 579]{DouVers}), we have universal morphisms
	\begin{equation} \label{eq:univ chart 2}
	\begin{split}
	(\underline{\beta}_{i} &: P_i \rightarrow \mathcal{O}_{p^{*}\underline{Y}^{\circ}_{i}})_{i\in I}\\
	(\underline{\varphi}^{i}_{j}:=(\phi^{i}_{j,0},\underline{\eta}^{i}_j)&: P_{i}\to P_{j}\oplus \O^{\times}_{p^{*}\underline{Y}^{\circ}_{j}})_{j\in I_1\cup I_2,i\in\partial j}
	\end{split}.
	\end{equation}
	By construction, they satisfy
	\[\big((\underline{\beta}_{j}\cdot\Id_{\O^{\times}_{\underline{Y}^{\circ}_{j}}})\circ\underline{\varphi}^{i}_{j}= \underline{\beta}_i|_{\underline{Y}^{\circ}_{j}}\big)_{j\in I_1\cup I_2, i\in\partial j}.\]
	On the other hand, we have that the space $\mathfrak{X}$ is canonically isomorphic to
	\[\coprod_{i\in I_0}\underline{Y}^{'}_{i}/\mathcal{R},\]
	where $\mathcal{R}(x,x^{'})$ if $x\in \underline{Y}'_{i}$ and  $x'\in \underline{Y}'_{i'}$ are such that there exists $j\in I_1$ and $y\in \underline{Y}'_{j}$ with $dj=(i,i')$, $\underline{g}^{j}_{i}(y)=x$ and $\underline{g}^{j}_{i'}(y)=x'$ (see \cite[p. 592]{DouVers}). Therefore, $\mathfrak{X}_{\log}$ is canonically isomorphic to
	\[\coprod_{i\in I_0}p^{*}\underline{Y}^{'}_{i}/\mathcal{R}.\]
	Hence, the collection of universal morphisms $((\underline{\beta}_{i}),(\underline{\varphi}^{i}_{j}))$ defines a pre-log atlas (see Definition $\ref{def:pre-log atlas}$) on $\mathfrak{X}_{\log}$, which glues to a fine log structure $\mathcal{M}_{\mathfrak{X}_{\log}}$ on $\mathfrak{X}_{\log}$ (see Subsection $\ref{sub:gluing charts}$). 
	\endproof
	We show that the universal family of log puzzles $(\mathfrak{X}_{\log},\mathcal{M}_{\mathfrak{X}_{\log}})\to \mathfrak{Z}^{\log}$ gives a \textit{complete} deformation of $(X_0,\mathcal{M}_{X_0})$.
	To do that, we introduce the notion of \textit{log cuirasse}. We recall that if $S$ is a Banach analytic space, $X$ a Banach analytic space proper and anaflat over $S$ and $q$ a relative cuirasse on $X$, then we get a morphism $\varphi_q:S\to\mathfrak{Z}$ ($\ref{def:morph associated to a cuirasse}$), a Banach analytic space $X_{\varphi_q}$ over $S$ obtained by gluing the pieces of the puzzle $z_q$ associated to $q$ ($\ref{prop:puzzle associated}$), and an $S$-isomporphism $\alpha_q:X_{\varphi_q}\to X$ ($\ref{prop:isomorphism alpha}$).
	Now, let $(X_0,\mathcal{M}_{X_0})$ be a compact fine log complex space admitting a collection of \text{directed} log charts $((\theta_i:P_i\to\mathcal{M}_{U_i})_{i\in I},(\varphi^{i}_j:=(\phi^{i}_j,\eta^i_j):P_i\to P_j\oplus \mathcal{O}^{\times}_{U_j})_{j\in I_1\cup I_2, i\in \partial j})$ (see Definition $\ref{def:compatible log charts}$).
	We assume that $\phi^{i}_j$ coincide with the $\phi^{i}_{j,0}$ given by Definition $\ref{def:discrete part}$.
	Let $q_0\in\mathcal{Q}(X_0)$ be a  cuirasse on $X_0$. We have an isomorphism ($\ref{prop:isomorphism alpha}$)
	\[\alpha_{q_0}:X_{\varphi_{q_0}}\to X_0.\]
	\begin{definition}\label{def:pull back log space}
		We naturally get a fine log structure on $X_{\varphi_{q_0}}$ via \[\mathcal{M}_{X_{\varphi_{q_0}}}:=\alpha_{q_0}^{*}\mathcal{M}_{X_0}.\]
	\end{definition}
	\begin{definition}\label{def:log cuirasse}
		A \textit{log cuirasse} $q_0^{\loga}$ on $(X_0, \mathcal{M}_{X_0})$ is a pair given by a cuirasse $q_0=(Y_i,f_i)_{i\in I}$ on $X_0$ and a  collection of directed log charts $((\theta_i:P_i\to\mathcal{M}_{X_{\varphi_{q_0}}}|_{Y^{\circ}_i}), (\eta^{i}_j:P_i\to O^{\times}_{Y^{\circ}_j}))$ on $(X_{\varphi_{q_0}},\mathcal{M}_{X_{\varphi_{q_0}}})$ (Definition $\ref{def:compatible log charts}$). We denote the set of log cuirasses on $(X_0, \mathcal{M}_{X_0})$ by $\mathcal{Q}(X_0,\mathcal{M}_{X_0})$.
		\begin{remark}
			In Definition $\ref{def:log cuirasse}$ we need to give the set of comparison morphisms $(\eta^{i}_{j})$ in order to define, in Definition $\ref{def:log puzzle associated}$, the log puzzle \textit{associated} to a log cuirasse.
		\end{remark}

	\end{definition}		
	
	Analogously to the classical case ($\ref{rel cuirasses}$), we can define the notion of relative log cuirasse. Let $S$ be a Banach analytic space and $(X,\mathcal{M}_{X})$ a fine log Banach analytic space proper and anaflat over $S$. 
	Given the local nature of the problem, we can assume that $(X,\mathcal{M}_{X})$ can be covered by finitely many log charts $(\theta_{i}:P_i\rightarrow \mathcal{M}_{U_i})_{i\in I}$ such that $U_{i}\cap X(s)\neq\emptyset$, for each $i\in I$ and $s\in S$. 
	
	\begin{definition}
		Let $S$ be a Banach analytic space and $(X,\mathcal{M}_{X})$ a fine log Banach analytic space proper and anaflat over $S$. We define the set of relative log cuirasses on $(X,\mathcal{M}_{X})$ over $S$ by
		\[\mathcal{Q}_{S}(X,\mathcal{M}_{X}):=\{(s,q)\lvert s\in S, q\in \mathcal{Q}(X(s),\mathcal{M}_{X}\arrowvert_{X(s)})\}= \coprod_{s\in S} \mathcal{Q}(X(s),\mathcal{M}_{X}\arrowvert_{X(s)}) .\]
		
	\end{definition}
	
	\begin{definition}\label{def:relative log cuirasse}
		We call a section $q^{\loga}:S\rightarrow \mathcal{Q}_{S}(X,\mathcal{M}_{X})$, of the canonical projection $\pi:\mathcal{Q}_{S}(X,\mathcal{M}_{X})\to S$,  a \textit{relative log cuirasse} on $(X,\mathcal{M}_{X})$ over $S$.  
	\end{definition}
	
	\begin{definition}\label{def:log triangularly priviliged}
		Let $(X_0,\mathcal{M}_{X_0})$ be a compact fine log complex space. A log cuirasse $q_0^{\loga}$ on $(X_0,\mathcal{M}_{X_0})$ is called \textit{triangularly privileged} if the underlying cuirasse $q_0\in\mathcal{Q}(\mathfrak{I};X_0)$ on $X_0$ is triangularly privileged (\ref{triang}). 
	\end{definition}
	
	Since every compact complex space $X_0$ admits a triangularly privileged cuirasse,  every compact fine log complex space $(X_0,\mathcal{M}_{X_0})$ admits a triangularly privileged log cuirasse. The set of log cuirasses can be endowed with the structure of a Banach analytic space in a neighborhood of a triangularly privileged log cuirasse. To prove it, we need the following three Lemmas.
	\begin{lemma}\label{lem:lemma 1 Banach}
		Let $(X,\mathcal{M}_X)$ be a fine log Banach analytic space over a Banach analytic space $S$. Let $q_0=(Y_{i,0},f_{i,0})$ be a triangularly privileged cuirasse on the central fibre $(X_0,\mathcal{M}_{X_0})$ over $s_0\in S$. Assume that $\Gamma(Y_{i,0},f^{-1}_{i,0}\overline{\mathcal{M}}_{X_0})$ is globally generated, for each $i\in I$. Then there exists a local relative cuirasse $q=(Y_i,f_i)$ on $X$ defined on a neighborhood $S'$ of $s_0$ in $S$, such that $\Gamma(Y_i,f^{-1}_i\overline{\mathcal{M}}_{X})$ is globally generated, for each $i\in I$.
	\end{lemma}
	\proof
	Since $q_0$ is triangularly privileged, there exists a local relative cuirasse $q=(Y_i,f_i)$ on $X$ defined in a neighborhood $S'$ of $s_0$ in $S$ (\cite[p. 585, Proposition 2]{DouVers}). Now, up to shrinking $Y_{i,0}$, for each $i\in I$ and for each $y\in Y_{i,0}$ there exists an open set $U_y=V_y\times W_y\subset S\times K_i$, $y\in U_y$, such that the canonical map $\Gamma(U_y,f_i^{-1}\overline{\mathcal{M}}_X)\to\overline{\mathcal{M}}_{X,f_i(y)}$ is an isomorphism. Since each $Y_i$ is compact, we can find a finite set $J$  and finitely many subsets $U_j$ such that $\Gamma(\bigcup_{j\in J} U_j,f_i^{-1}\overline{\mathcal{M}}_X)$ is isomorphic to $\Gamma(Y_i,f_i^{-1}\overline{\mathcal{M}}_X)$. Hence, possibly after shrinking $S'$, we can assume that for each $j\in J$, $V_j=S'$ and $f^{-1}_i(X)\subset\bigcup_{j\in J}U_j$. Thus, we get that for each $i\in I$, $\Gamma(Y_i,f^{-1}_i\overline{\mathcal{M}}_{X})$ is isomorphic to $\Gamma(Y_{i,0},f^{-1}_{i,0}\overline{\mathcal{M}}_{X})$.
	\endproof
	A holomorphic line bundle with $c_1=0$ is topologically trivial, hence analytically isomorphic to the trivial line bundle by the following Lemma \ref{lem:Lemma 2 Banach}.
	\begin{lemma}(\cite[p. 268]{GraLine})\label{lem:Lemma 2 Banach}
		Let $X$ be a Stein space and $\mathcal{L}$ a $\mathcal{O}^{\times}_X$-torsor. If $c_1(\mathcal{L})=0$, then $\mathcal{L}$ is trivial.
	\end{lemma}

	
	\begin{lemma}\label{lem:Lemma 3 Banach}
		Let $(X,\mathcal{M}_X)$ be a fine log Banach analytic space. Assume that $P:=\Gamma(X,\overline{\mathcal{M}}_X)$ is globally generated and torsion free. Assume that for each $\overline{m}\in \Gamma(X,\overline{\mathcal{M}}_X)$, the torsor $\mathcal{L}_{\overline{m}}=\kappa^{-1}(\overline{m})$, with $\kappa:\mathcal{M}_X\to\overline{\mathcal{M}}_X$ the canonical map, is trivial. Then there exists a chart $P\to \Gamma(X,\mathcal{M}_X)$.
	\end{lemma}
	\proof
	Let $p_1,...,p_r \in P$ be generators, that is we have a surjective map $\mathbb{N}^r\to \Gamma(X,\overline{\mathcal{M}}_X)$ sending $e_i$ to $p_i$. For each $i\in\{1,...,r\}$, choose a section $m_i\in\mathcal{L}_{\overline{m}}$.  We obtain a chart $\phi:\mathbb{N}^r\to \Gamma(X,\mathcal{M}_X)$.  Now, we want to modify $\phi$ so that it factors through $P$. Let $K:=\ker(\mathbb{Z}^r\to P^{gp})$, we have $P=\mathbb{N}^r/K$. We get the following exact sequence
	\[0\to K\to \mathbb{Z}^r\to P^{gp}\to 0.\]
	Since, by assumption, $P$ is torsion free, we can find a section $\pi:\mathbb{Z}^r\to K$. Set $h_i:=\phi^{gp}(\pi(p_i))$. Clearly, if $\sum a_ip_i=\sum b_jp_j$, for $a_i, b_j
	\geq 0$, then it holds $\prod h^{a_i}_i=\prod h^{b_j}_j$ in $\Gamma(X,\mathcal{M}^{gp}_X)$. We get a chart by
	
	\begin{equation}
	\begin{split}
	\tilde{\psi}:\mathbb{N}^r&\to\Gamma(X,\mathcal{M}_X)\\
	e_i&\mapsto h_i^{-1}m_i\\ 
	\end{split}.
	\end{equation}
	
	Now, let $\tilde{\psi}^{gp}:\mathbb{Z}^r\to\Gamma(X,\mathcal{M}^{gp}_X)$ and  $\sum a_ie_i\in K$. If $\tilde{ \psi}^{gp}(\sum a_ie_i)=1$, we get that $\tilde{\psi}^{gp}$ induces a chart $\psi:P\to\Gamma(X,\mathcal{M}_X)$.
	Hence, assume $\sum a_ie_i\in K$, $a_i\in \mathbb{Z}$. Then $\tilde{ \psi}^{gp}(\sum a_ie_i)=\prod \tilde{ \psi}^{gp}(e_i)^{a_i}=\prod h^{-a_i}_im_i^{a_i}=\prod \phi^{gp}(\pi(p_i))^{-a_i}\phi^{gp}(e_i)^{a_i}=\phi^{gp}(\pi(-\sum a_ip_i)+\sum a_ie_i))=1$.

	\begin{proposition}\label{prop:log fibrewise cuirasse}
		Let $(X,\mathcal{M}_X)$ be a fine log Banach analytic space over a Banach analytic space $S$. Let $s_0\in S$ and $q^{\loga}_0$ a triangularly privileged log cuirasse on $(X(s_0),\mathcal{M}_{X(s_0)})$.
		Then the set of log cuirasses $\mathcal{Q}_S(X,\mathcal{M}_X)$ on $(X,\mathcal{M}_X)$ over $S$ can be endowed with the structure of a Banach analytic space in a neighborhood of $(s_0,q_0^{\loga})$.
	\end{proposition}
	\proof
	Let us consider the projection $\pi:\mathcal{Q}_S(X,\mathcal{M}_X)\to Q_S(X)$.  By Lemma $\ref{lem:lemma 1 Banach}$ and Lemma $\ref{lem:Lemma 2 Banach}$, we can use Lemma $\ref{lem:Lemma 3 Banach}$ and get the existence around $(s_0,q_0)$ of a local section $\rho:Q_S(X)\to \mathcal{Q}_S(X,\mathcal{M}_X)$, such that $\rho(s_0,q_0)=q^{\loga}_0$. Now, let $(s,q)\in Q_S(X)$, in a small neighborhood of $(s_0,q_0)$, and consider $\rho(s,q)\in \mathcal{Q}_S(X,\mathcal{M}_X)$.  We have that $\rho(s,q)=(s,q=(Y_i,f_i),(\theta_i),(\eta^{i}_j))$, where $(\theta_i),(\eta^{i}_j)$ is a directed collection of log charts on $(X_{\varphi_q},\mathcal{M}_{X_{\varphi_q}})$ (see Definition $\ref{def:log cuirasse}$). 
	Any other directed set of log charts $((\theta'_i),(\eta^{'i}_j))$ on $(X_{\varphi_q},\mathcal{M}_{X_{\varphi_q}})$ is obtained by $\theta'_i=\chi_i\cdot\theta_i$ and $\eta^{'i}_j=\chi_i^{-1}\cdot\chi_j\cdot\eta^i_j$, for morphisms  $\chi_i:P_i\to\mathcal{O}^{\times}_{Y^{\circ}_i}$, for $i\in I$. Therefore,  let $\underline{Y}_i$ be the universal $\mathcal{G}(K_i)$-anaflat subspace of $\mathcal{G}(K_i)\times K_i$, for $i\in I$ (\cite[p. 579]{DouVers}, \cite[pp. 258--259]{Lepot} and \cite[p. 183, Theorem 4.13]{Pourpriv}). We can define a map
	\[\gamma:\mathcal{Q}_S(X)\times_{\prod \mathcal{G}(K_i)\times S}\fprods_{i\in I} \mathfrak{Mor}_{\mathcal{G}(K_i)\times S}(\underline{Y}_i\times S, \Spec\mathbb{C}[P^{gp}_{i}]\times \mathcal{G}(K_i)\times S)\to \mathcal{Q}_S(X,\mathcal{M}_X)\]
	via
	\[(s,q,(\chi_i))\mapsto(s,q,(\chi_i\cdot\theta_i), (\chi_i^{-1}\cdot\chi_j\cdot\eta^i_j)),\]
	which defines a structure of Banach analytic space on $\mathcal{Q}_S(X,\mathcal{M}_X)$ in a neighborhood of $(s_0,q^{\loga}_0)$.
	
	\endproof
	\begin{proposition}\label{prop:log triangularly smoothness}
		Let $(X,\mathcal{M}_X)$ be a fine log Banach analytic space proper and anaflat over a Banach analytic space $S$. Let $s\in S$ and $q^{\loga}(s)$ a triangularly privileged log cuirasse on $(X(s),\mathcal{M}_{X(s)})$. Then
		\[\pi:\mathcal{Q}_S(X,\mathcal{M}_X)\to S\]
		is smooth in a neighborhood of $q^{\loga}(s)$.
	\end{proposition}
	\proof
	Let $q^{\loga}(s)$ be a triangularly privileged log cuirasse on $(X(s),\mathcal{M}_{X(s)})$. By Proposition $\ref{prop:log fibrewise cuirasse}$, we have that in a neighborhood of  $(s,q^{\loga}(s))$, the space $\mathcal{Q}_{S}(X,\mathcal{M}_{X})$ is isomorphic to 
	\[\mathcal{Q}_S(X)\times_{\prod \mathcal{G}(K_i)\times S}\fprods_{i\in I} \mathfrak{Mor}_{\mathcal{G}(K_i)\times S}(\underline{Y}_i\times S, \Spec\mathbb{C}[P^{gp}_{i}]\times \mathcal{G}(K_i)\times S).\]
	Let $q(s)=(Y_i,f_i)$	be the triangularly privileged cuirasse on $X(s)$ underlying $q^{\loga}(s)$. By \cite[p. 589, Corollary 2]{DouVers}, $\pi:\mathcal{Q}_{S}(X)\to S$ is smooth in a neighborhood of $(s,q(s))$. Furthermore, by \cite[p. 585, Proposition 2]{DouVers}, we have that \[\mathfrak{Mor}_{\mathcal{G}(K_i)\times S}(\underline{Y}_i\times S, \Spec\mathbb{C}[P^{gp}_{i}]\times \mathcal{G}(K_i)\times S)\to S\]
	is smooth in a neighborhood of $(s,Y_i,\chi_i)$. Hence, the statement follows.
	\endproof

	Analogously to the classical case (\ref{prop:puzzle associated}), we define the notion of log puzzle \textit{associated} to a log cuirasse.
	Let $(X_0,\mathcal{M}_{X_0})$ be a compact fine log complex space. Let $q_0$ be a cuirasse on $X_0$. By Definition $\ref{def:pull back log space}$, we get a compact fine log complex space $(X_{\varphi_{q_0}},\mathcal{M}_{X_{\varphi_{q_0}}})$, which is isomorphic to $(X_0,\mathcal{M}_{X_0})$. Let $q^{\loga}_0=(q_0,(\theta_i),(\eta^{i}_j))$ be a log cuirasse on $(X_0,\mathcal{M}_{X_0})$ (see Definition $\ref{def:log cuirasse}$). Let $\alpha_{X_{\varphi_{q_0}}}:\mathcal{M}_{X_{\varphi_{q_0}}}\to\mathcal{O}_{X_{\varphi_{q_0}}}$ be the structure log morphism  and $z_{q_0}\in\mathfrak{Z}$ the puzzle associated to $q_0$.

	\begin{definition}\label{def:log puzzle associated}
		We call \[z_{q^{\loga}_0}:=(z_{q_0},(\alpha_{X_{\varphi_{q_0}}}\circ\theta_i),(\eta^{i}_j
		))\] the log puzzle \textit{associated} to $q^{\loga}_0$.
	\end{definition}
	Clearly, $z_{q^{\loga}_0}\in\mathfrak{Z}^{\log}$ (Definition $\ref{def:log puzzle}$).	
	Let $(X,\mathcal{M}_{X})$ be a fine log Banach analytic space proper and anaflat over a Banach analytic space $S$. Let $q^{\loga}$ be a relative log cuirasse on $(X,\mathcal{M}_{X})$ over $S$. 
	\begin{definition}\label{def:log morph associated}
		We can define a morphism
		\begin{equation}
		\begin{split}
		\varphi_{q^{\loga}}:S&\to\mathfrak{Z}^{\log}\\
		s&\mapsto z_{q^{\loga}(s)}.	\\
		\end{split} 
		\end{equation}

	\end{definition}
	
	%
	
	
	Let $q^{\loga}=(q,(\theta_i),(\eta^i_j))$ be a log cuirasse on $(X,\mathcal{M}_{X})$ over $S$.
	Let $(X_{\varphi_q},\mathcal{M}_{X_{\varphi_q}})$ given by  Definition $\ref{def:pull back log space}$ and $\alpha_{X_{\varphi_q}}:\mathcal{M}_{X_{\varphi_q}}\to\mathcal{O}_{X_{\varphi_q}}$ the structure log morphism. For each $i\in I$, let  $\mathcal{M}^{a}_{X_{\varphi_q},i}$ be the log structure associated to the pre-log structure $\alpha_{X_{\varphi_q}}\circ\theta_i$.  The collection of log structures $(\mathcal{M}_{X_{\varphi_q},i}^{a})$ glues to a log structure $\mathcal{M}_{X_{\varphi_q}}^{a}$ on $X_{\varphi_q}$ (Subsection \ref{sub:gluing charts}).
	
	
	\begin{definition}
		We set\[(X_{\varphi_{q^{\loga}}},\mathcal{M}_{X_{\varphi_{q^{\loga}}}}):=(X_{\varphi_{q}},\mathcal{M}^{a}_{X_{\varphi_q}}).\]
	\end{definition}

	The fine log Banach analytic space $(X_{\varphi_{q^{\loga}}},\mathcal{M}_{X_{\varphi_{q^{\loga}}}})$ is obtained by gluing the pieces of the log puzzle $z_{q^{\loga}}$ associated to the cuirasse $q^{\loga}$.
	\begin{proposition}\label{prop:log iso alpha}
		Let $(X,\mathcal{M}_X)$ be a fine log Banach analytic space proper and anaflat over a Banach analytic space $S$.  Let $q^{\loga}$ be a relative log cuirasse on $(X,\mathcal{M}_X)$ over $S$. 
		Then, there exists a log $S$-isomorphism

		\[	\begin{tikzpicture}[descr/.style={fill=white,inner 		
			sep=2.5pt}]
		\matrix (m) [matrix of math nodes, row sep=3em,
		column sep=3em]
		{ (X_{\varphi_{q^{\loga}}},\mathcal{M}_{X_{\varphi_{q^{\loga}}}}) &  &(X,\mathcal{M}_X)\\
			& S & \\ };
		\path[->,font=\scriptsize]
		(m-1-1) edge node[above] {$\alpha_{q^{\loga}} $} (m-1-3) 
		(m-1-1) edge node[left] {$ $} (m-2-2)
		(m-1-3) edge node[below] {$ $} (m-2-2)
		;
		\end{tikzpicture}.\]

	\end{proposition}
	\proof
	By Proposition $\ref{prop:isomorphism alpha}$, we have an $S$-isomorphism $\alpha_q:X_{\varphi_q}\to X$. Moreover, we have $\mathcal{M}_{X_{\varphi_q}}:=\alpha_q^{*}\mathcal{M}_{X}$ (see Definition $\ref{def:log cuirasse}$). Hence, $\alpha_q$ induces a $S$-log isomorphism $\alpha_q:(X_{\varphi_q},\mathcal{M}_{X_{\varphi_q}})\to(X,\mathcal{M}_X)$. Now, the log cuirasse $q^{\loga}$ gives us a collection of directed log charts $((\theta_i),(\eta^i_j))$ for $\mathcal{M}_{X_{\varphi_q}}$. Let $\alpha_{X_{\varphi_q}}:\mathcal{M}_{X_{\varphi_q}}\to\mathcal{O}_{X_{\varphi_q}}$ be the structure log morphism and $\mathcal{M}^{a}_{X_{\varphi_q},i}$ the log structure associated to the pre-log structure $\alpha_{X_{\varphi_q}}\circ\theta_i$, for each $i\in I$. By the definition of log chart (\cite[p. 249]{Ogus18}), we have an isomorphism $\alpha^{\flat}_i:\mathcal{M}_{X_{\varphi_q},i}^{a}\to\mathcal{M}_{X_{\varphi_q}}$. Then the collection of log structures $(\mathcal{M}_{X_{\varphi_q},i}^{a})$, together with the isomorphisms $(\alpha^{\flat}_i)$, glues to a log structure $\mathcal{M}_{X_{\varphi_q}}^{a}$ on $X_{\varphi_q}$, together with an isomorphism $\alpha^{\flat}:\mathcal{M}_{X_{\varphi_q}}^{a}\to\mathcal{M}_{X_{\varphi_q}}$ (see Subsection  \ref{sub:gluing charts}). Hence, set $(X_{\varphi_{q^{\loga}}},\mathcal{M}_{X_{\varphi_{q^{\loga}}}}):=(X_{\varphi_{q}},\mathcal{M}^{a}_{X_{\varphi_q}})$ and $\alpha:=(\Id, \alpha^{\flat})$, we get an isomorphism $\alpha:(X_{\varphi_{q^{\loga}}},\mathcal{M}_{X_{\varphi_{q^{\loga}}}})\to (X_{\varphi_{q}},\mathcal{M}_{X_{\varphi_q}})$. Set $\alpha_{q^{\dagger}}:=\alpha_q\circ\alpha$.
	\endproof
	\begin{remark}
		Clearly,  \[(X_{\varphi_{q^{\loga}}},\mathcal{M}_{X_{\varphi_{q^{\loga}}}})=\varphi_{q^{\loga}}^{*}(\mathfrak{X}_{\log},\mathcal{M}_{\mathfrak{X}_{\log}}).\]
	\end{remark}
	We are ready to prove the existence of an infinite-dimensional complete deformation of a fine compact log complex space $(X_0,\mathcal{M}_{X_0})$. With the due modifications, the proof of Theorem $\ref{thm:inf dim log versal def}$ is identical to the proof of Theorem $\ref{thm:douady}$ (\cite[p. 592]{DouVers}).
	Let $q_0=(Y_{i,0},f_{i,0})$ be a triangularly privileged cuirasse on $X_0$ and $((\theta_{i,0}),(\eta^{i}_{j,0}))$ the collection of directed log charts on $(X_{0},\mathcal{M}_{X_0})$ as in Definition $\ref{def:discrete part}$. Then, \[q^{\loga}_{0}:=(q_0,(f^{*}_{i,0}\theta_i),(f^{*}_{j,0}\eta^{i}_j))\] is a triangularly privileged log cuirasse on $(X_{0},\mathcal{M}_{X_0})$ (see Definition $\ref{def:log triangularly priviliged}$). 
	Let $z_{q^{\loga}_{0}}$ be the log puzzle associated to $q^{\loga}_{0}$ (see Definition $\ref{def:log puzzle associated}$). Let $(\mathfrak{X}_{\log},\mathcal{M}_{\mathfrak{X}_{\log}})\to \mathfrak{Z}^{\log}$ be the universal space of log puzzles (see Proposition $\ref{prop:log moduli}$) and  \[\alpha_{q^{\loga}_{0}}:(\mathfrak{X}_{\log}(z_{q^{\loga}_{0}}),\mathcal{M}_{\mathfrak{X}_{\log}(z_{q^{\loga}_{0}})})\rightarrow (X_{0},\mathcal{M}_{X_0})\] the log isomorphism given by Proposition $\ref{prop:log iso alpha}$.
	
	
	\begin{theorem}\label{thm:inf dim log versal def}
		The triple $((\mathfrak{Z}^{\log}, z_{q^{\loga}_0}), (\mathfrak{X}_{\log},\mathcal{M}_{\mathfrak{X}_{\log}}), \alpha_{q^{\loga}_{0}})$ is a complete deformation of $(X_{0},\mathcal{M}_{X_0})$. 
	\end{theorem}
	\proof
	Let $(X,\mathcal{M}_{X})$ be a fine log Banach analytic space proper and anaflat over a Banach analytic space $S$. Let $s_0\in S$ and $i:(X(s_0),\mathcal{M}_{X(s_0)})\to (X_0,\mathcal{M}_{X_0})$ a log isomorphism. Since $i^{*}q^{\loga}_{0}$ is a triangularly privileged log cuirasse (see Definition $\ref{def:log triangularly priviliged}$) on $(X(s_{0}),\mathcal{M}_{X(s_0)})$, we have that $\mathcal{Q}_S(X,\mathcal{M}_X)$ is smooth over $S$ in a neighborhood of $i^{*}q^{\loga}_{0}$ (Proposition $\ref{prop:log triangularly smoothness}$). Therefore there exists a local relative log cuirasse $q^{\loga}$ on $(X,\mathcal{M}_X)$ defined in a neighborhood $S'$ of $s_{0}$ in $S$. 	
	Hence, taking $\varphi_{q^{\loga}}:S\to \mathfrak{Z}^{\log}$ (Definition $\ref{def:log morph associated}$) and the $S'$-isomorphism \[\alpha_{q^{\loga}}|_{S'}:(X_{\varphi_{q^{\loga}}},\mathcal{M}_{X_{\varphi_{q^{\loga}}}})\rightarrow (X,\mathcal{M}_X)\] (Proposition $\ref{prop:log iso alpha}$), the statement follows.
	\endproof

	\subsection{Finite dimensional reduction}
	\label{finite dimensional reduction}

		Let $(\mathfrak{X}_{\log},\mathcal{M}_{\mathfrak{X}_{\log}})\to \mathfrak{Z}^{\log}$ be the complete deformation of $(X_0,\mathcal{M}_{X_0})$ given by Theorem $\ref{thm:inf dim log versal def}$ and $\mathcal{Q}_{\mathfrak{Z}^{\log}}(\mathfrak{X}_{\log},\mathcal{M}_{\mathfrak{X}_{\log}})$ the space  of relative log cuirasses on $(\mathfrak{X}_{\log},\mathcal{M}_{\mathfrak{X}_{\log}})$ over $\mathfrak{Z}^{\log}$ (Definition $\ref{def:log cuirasse}$). Since the finite-dimensional reduction is performed on the Banach analytic space $\mathcal{Q}_{\mathfrak{Z}^{\log}}(\mathfrak{X}_{\log},\mathcal{M}_{\mathfrak{X}_{\log}})$, which does not come endowed with a non-trivial log structure, the finite-dimensional reduction in the log setting is identical to the one in the classical setting (see \cite[pp. 593--599]{DouVers} and \cite[pp. 20--46]{Stie}).  
	In what follows, we give an account of the main steps of the finite-dimensional reduction procedure (in the log setting). For more details, see \cite[pp. 90--100]{Raf}.
	\smallskip
	
	We recall from subection $\ref{Infinite dimensional construction}$ that the space $\mathfrak{Z}^{\log}$ (Definition $\ref{def:log puzzle set}$) parametrizes all log puzzles $z^{\loga}$ of type $\mathfrak{I}$ (Definition $\ref{def:log puzzle}$). Each fibre $(\mathfrak{X}_{\log},\mathcal{M}_{\mathfrak{X}_{\log}})(z^{\loga})$ of the map  $(\mathfrak{X}_{\log},\mathcal{M}_{\mathfrak{X}_{\log}})\to \mathfrak{Z}^{\log}$ 
	is obtained by gluing the \enquote{pieces}, $(Y_i)_{i\in I}$ and  $(\beta_i:P_i\rightarrow\mathcal{O}_{Y^{\circ}_i})_{i\in I}$, of the log puzzle $z^{\loga}$. 
	
	Each point in $\mathcal{Q}_{\mathfrak{Z}^{\log}}(\mathfrak{X}_{\log},\mathcal{M}_{\mathfrak{X}_{\log}})$ is a pair $(z^{\loga},q^{\loga})$, where $z^{\loga}\in\mathfrak{Z}^{\log}$ is a log puzzle and $q^{\loga}$ is a log cuirasse on the fibre $(\mathfrak{X}_{\log},\mathcal{M}_{\mathfrak{X}_{\log}})(z^{\loga})$. To  the log cuirasse $q^{\loga}$  we can naturally associate another log puzzle $z_{q^{\loga}}\in\mathfrak{Z}^{\loga}$ (Definition $\ref{def:log puzzle associated}$).
	In principle,  $z^{\loga}\neq z_{q^{\loga}}$ although \[(\mathfrak{X}_{\log},\mathcal{M}_{\mathfrak{X}_{\log}})(z^{\loga})\simeq(\mathfrak{X}_{\log},\mathcal{M}_{\mathfrak{X}_{\log}})(z_{q^{\loga}}),\] 
	(Proposition $\ref{prop:log iso alpha}$). However, we can consider the subspace $Z^{\log}\subset\mathcal{Q}_{\mathfrak{Z^{\log}}}(\mathfrak{X}_{\log},\mathcal{M}_{\mathfrak{X}_{\log}})$ defined by selecting, in each fibre $\mathcal{Q}((\mathfrak{X}_{\log},\mathcal{M}_{\mathfrak{X}_{\log}})(z^{\loga}))$ of the canonical projection $\pi:\mathcal{Q}_{\mathfrak{Z}^{\log}}(\mathfrak{X}_{\log},\mathcal{M}_{\mathfrak{X}_{\log}})\to\mathfrak{Z}^{\log}$,  all log cuirasses $q^{\loga}$ on $(\mathfrak{X}_{\log},\mathcal{M}_{\mathfrak{X}_{\log}})(z^{\loga})$ whose associated log puzzle $z_{q^{\loga}}$ coincides exactly with $z^{\loga}$. 
	More precisely, there exists a canonical relative log cuirasse $\mathfrak{q}^{\loga}$  on
	\begin{equation}\label{log can}
	\pi^{*}(\mathfrak{X}_{\log},\mathcal{M}_{\mathfrak{X}_{\log}})\to \mathcal{Q}_{\mathfrak{Z}^{\log}}(\mathfrak{X}_{\log},\mathcal{M}_{\mathfrak{X}_{\log}}),
	\end{equation} 
	see \cite[p. 593]{DouVers}, \cite[p.267]{Lepot} and \cite[p. 90]{Raf}. 	By Definition $\ref{def:log morph associated}$, we get an associated morphism
	\begin{equation}\label{log can cuirasse}
	\begin{split}
	\varphi_{\mathfrak{q}^{\loga}}:\mathcal{Q}_{\mathfrak{Z}^{\log}}(\mathfrak{X}_{\log},\mathcal{M}_{\mathfrak{X}_{\log}})&\to\mathfrak{Z}^{\log}\\
	(z^{\loga},q^{\loga})&\mapsto z_{q^{\loga}}\\
	\end{split}.
	\end{equation}
	Then, the subspace $Z^{\log}$ is obtained as the kernel of the double arrow $(\pi,\varphi_{\mathfrak{q}^{\loga}})$:
	\begin{equation}\label{eq:log space Z def}
	Z^{\log}:=\ker(\pi,\varphi_{\mathfrak{q}^{\loga}})\subset\mathcal{Q}_{\mathfrak{Z}^{\log}}(\mathfrak{X}_{\log},\mathcal{M}_{\mathfrak{X}_{\log}}).
	\end{equation}
	The space $Z^{\log}$ parametrizes all log cuirasses on compact fine log complex spaces \enquote{close} to $(X_0,\mathcal{M}_{X_0})$. This space is not as pathological as $\mathfrak{Z^{\log}}$ (see \cite[p. 590, Remark]{DouVers}) and it still gives a complete deformation of $(X_0,\mathcal{M}_{X_0})$.
	Indeed, given any log Banach analytic space $(X,\mathcal{M}_X)$ proper and anaflat over a Banach analytic space $S$, we get a map from the space of relative log cuirasses $\mathcal{Q}_S(X,\mathcal{M}_X)$ into the space of log puzzles $\mathfrak{Z}^{\log}$ by
	\begin{equation}
	\begin{split}
	\varphi_{(X,\mathcal{M}_X)/S}:\mathcal{Q}_S(X,\mathcal{M}_X)&\to\mathfrak{Z}^{\log}\\
	(s,q^{\loga})&\mapsto z_{q^{\loga}}.\\
	\end{split}
	\end{equation}
	If $\sigma^{\loga}:S\to \mathcal{Q}_S(X,\mathcal{M}_X)$ is a relative log cuirasse on $(X,\mathcal{M}_X)$ over $S$, that is a section of the projection $\pi:\mathcal{Q}_S(X,\mathcal{M}_X)\to S$, the composition \[\varphi_{\sigma^{\loga}}:=\varphi_{(X,\mathcal{M_X})/S}\circ \sigma^{\loga}:S\to\mathfrak{Z}^{\log}\] is a morphism satisfying the completeness property (see Definition $\ref{def:log morph associated}$ and Theorem $\ref{thm:inf dim log versal def}$).
	Indeed, for each $s\in S$, the fibre $(\mathfrak{X}_{\log},\mathcal{M}_{\mathfrak{X}_{\log}})(\varphi_{\sigma^{\loga}}(s))$ is isomorphic to the fibre $(X,\mathcal{M}_X)(s)$ via an isomorphism $\alpha_{\sigma^{\loga}}$ (Proposition $\ref{prop:log iso alpha}$). Identifying these two isomorphic fibres, we get a map
	\begin{equation}\label{def:map psi simple}
	\begin{split}
	\psi_{\sigma^{\loga}}:S&\to Q_{\mathfrak{Z}^{\log}}(\mathfrak{X}_{\log},\mathcal{M}_{\mathfrak{X}_{\log}})\\
	s&\mapsto (z_{\sigma^{\loga}(s)},\sigma^{\loga}(s))\\
	\end{split}.
	\end{equation}
	For more details, see \cite[pp. 70--71 and pp. 90--91]{Raf}. We can draw the following commutative diagram: 
	\[	\begin{tikzpicture}[descr/.style={fill=white,inner 		
		sep=2.5pt}]
	\matrix (m) [matrix of math nodes, row sep=3em,
	column sep=3em]
	{ & & Q_{\mathfrak{Z}^{\log}}(\mathfrak{X}_{\log},\mathcal{M}_{\mathfrak{X}_{\log}})\\
		S& \mathcal{Q}_S(X,\mathcal{M}_X) & \mathfrak{Z}^{\log} \\};
	\path[->,font=\scriptsize]
	(m-2-1) edge node[below] {$\sigma^{\loga}$} (m-2-2) 
	(m-2-2) edge node[below] {$\varphi_{(X,\mathcal{M}_X)/S}$} (m-2-3) 
	(m-2-1) edge node[above] {$\psi_{\sigma^{\loga}}$} (m-1-3)
	(m-1-3) edge node[right] {$\pi$} (m-2-3);
	\end{tikzpicture}.\]
	In fact, $\psi_{\sigma^{\loga}}$ is the unique morphism from $S$ to $Q_{\mathfrak{Z}^{\log}}(\mathfrak{X}_{\log},\mathcal{M}_{\mathfrak{X}_{\log}})$ making the above diagram commutative 
	(see \cite[p. 593]{DouVers}). 
	By construction, $\psi_{\sigma^{\loga}}$ factors through $Z^{\log}\subset Q_{\mathfrak{Z}^{\log}}(\mathfrak{X}_{\log},\mathcal{M}_{\mathfrak{X}_{\log}})$ and it is used to prove Proposition \ref{prop:versal Z}.	 
	Let
	\[i:Z^{\log}\hookrightarrow\mathcal{Q}_{\mathfrak{Z}^{\log}}(\mathfrak{X}_{\log},\mathcal{M}_{\mathfrak{X}_{\log}})\] be the canonical injection, we set \[(\mathfrak{X}_{Z^{\log}},\mathcal{M}_{\mathfrak{X}_{Z^{\log}}}):=i^{*}\pi^{*}(\mathfrak{X}_{\log},\mathcal{M}_{\mathfrak{X}_{\log}}).\] 
	Let $q^{\loga}_0$ be a triangularly privileged log cuirasse on $(X_0,\mathcal{M}_{X_0})$ (see Definition $\ref{def:log triangularly priviliged}$) and $z_{q^{\loga}_0}\in\mathfrak{Z}^{\log}$ the associated log puzzle (Definition $\ref{def:log puzzle associated}$). We get a point  $(z_{q^{\loga}_0},q^{\loga}_0)$ in $Z^{\log}$.
	
	\begin{proposition}(\cite[p. 593, Theorem 1]{DouVers}, \cite[p. 31, Satz 1.15]{Stie}, \cite[p. 91, Proposition 5.26]{Raf})\label{prop:versal Z}
		The morphism $(\mathfrak{X}_{Z^{\log}},\mathcal{M}_{\mathfrak{X}_{Z^{\log}}})\to Z^{log}$ is a complete deformation of the compact fine log complex space $(X_0,\mathcal{M}_{X_0})$, in a neighborhood of $(z_{q^{\loga}_0},q^{\loga}_0)$.
	\end{proposition}
	In what follows,  we are going to decompose $Z^{\log}$ into a product $\Sigma^{\log}\times R^{\log}$, where $\Sigma^{\log}$ is a Banach manifold  and $R^{\log}$ is a finite dimensional complex analytic space, which will be our finite dimensional semi-universal deformation space.
	\medskip

	To do that, let us start by introducing the notion of \textit{extendable log cuirasse} (\cite[pp. 95--96]{Raf}) by adapting, to the log context,  Douady's notion of \textit{extendable cuirasse} (\cite[p. 594]{DouVers}). This is a fundamental tool to achieve finite dimensionality.

	\begin{definition}(\cite[p. 594]{DouVers})\label{def:extendable type}
		Let us consider two types of cuirasses, namely  $\mathfrak{I}=(I_{\bullet}, (K_{i}),(\tilde{K}_{i}),(K'_{i}))$ and $\hat{\mathfrak{I}}=(I_{\bullet}, (\hat{K}_{i}),(\hat{\tilde{K}}_{i}),(\hat{K'}_{i}))$  ($\ref{type cuirasse}$), which have the same underlying simplicial set. We write $\mathfrak{I}\Subset \hat{\mathfrak{I}}$, if
		$K_i\Subset \hat{K}_i$, $\tilde{K}_{i}\subset\hat{\tilde{K}}_{i}$ and $K'_{i}\subset \hat{K'}_{i}$. 
	\end{definition}
	
	Let $\hat{\mathfrak{I}}$ be a type of cuirasse and $\hat{q}^{\loga}$ a relative log cuirasse of type  $\hat{\mathfrak{I}}$ on a log Banach analytic space $(X,\mathcal{M}_X)$ proper and anaflat over $S$. Then, by slightly shrinking each polycylinder $\hat{K}_i$, $\hat{\tilde{K}}_{i}$ and $\hat{K}'_{i}$ in $\hat{\mathfrak{I}}$, we can get polycylinders $K_i$, $\tilde{K}_{i}$ and $K'_{i}$ respectively and hence a type of cuirasse $\mathfrak{I}$, such that $\mathfrak{I}\Subset\mathfrak{\hat{I}}$. 
	Then, \[q^{\loga}:=\hat{q}^{\loga}|_{\mathfrak{I}}\] is an \textit{extendable} relative log  cuirasse on  $(X,\mathcal{M}_X)$ over $S$. 
	
	If $\mathfrak{I}\Subset\hat{\mathfrak{I}}$ are two types of cuirasses, then we can construct the spaces of log puzzles $\mathfrak{Z}^{\log}$ and $\hat{\mathfrak{Z}}^{\log}$ of type   $\mathfrak{I}$ and $\hat{\mathfrak{I}}$ respectively (see Definition $\ref{def:log puzzle set}$). It can be shown (see \cite[p. 595]{DouVers} and \cite[p.44]{Stie}), that the restriction morphism
	\begin{equation}\label{restriction}
	j^{\loga}:\hat{\mathfrak{Z}}^{\log}\rightarrow \mathfrak{Z}^{\log}
	\end{equation}
	is compact (in the sense of \cite[p. 28]{DouThesis}). 
	This fact, together with the finite dimensionality results \cite[p. 29, Proposition 3]{DouThesis} and \cite[p. 271]{Lepot} (see, also, \cite[pp. 43-44]{Raf}), is used to prove Proposition \ref{prop:log rel finite}.
	\smallskip
	
	Set \[Q^{\log}_0:=\mathcal{Q}(X_0,\mathcal{M}_{X_0}),\]
	the space of log cuirasses on $(X_0,\mathcal{M}_{X_0})$ (see Definition $\ref{def:log cuirasse}$).
	By Proposition $\ref{prop:log triangularly smoothness}$ the projection $\pi:\mathcal{Q}_{\mathfrak{Z}^{\log}}(\X_{\log},\mathcal{M}_{\X_{\log}})\to \mathfrak{Z}^{\log}$ is smooth in a neighborhood of $(z_{q^{\loga}_0},q^{\loga}_0)$,
	hence we can opportunely choose (see \cite[p. 595]{DouVers}, \cite[p. 269]{Lepot} and \cite[p. 96]{Raf})  a local trivialization 
	\begin{equation}\label{trivial log}
	(\pi,\rho^{\loga}):\mathcal{Q}_{\mathfrak{Z}^{\log}}(\X_{\log},\mathcal{M}_{\X_{\log}})\to \mathfrak{Z}^{\log}\times Q^{\log}_0.
	\end{equation}
	Set
	\begin{equation}\label{eq: morph p log}
	p^{\loga}:=\rho^{\loga}|_{Z^{\log}}.
	\end{equation}
	Proposition $\ref{prop:log rel finite}$ is the log version of \cite[p. 596, Proposition 4]{DouVers}, \cite[p.43, Satz 1.33]{Stie} and \cite[p. 269, Lemma 1]{Lepot}. The same proof applies here likewise. For further details, see \cite[p. 96, Proposition 5.30]{Raf}.
	\begin{proposition}\label{prop:log rel finite}
		The morphism $p^{\loga}:Z^{\log}\to Q^{\log}_0$ is of relative finite dimension in a neighborhood of $(z_{q^{\loga}_0},q^{\loga}_0)$.
	\end{proposition}
	Thus, we get the existence of an embedding $\iota^{\loga}:Z^{\log}\hookrightarrow Q^{\log}_0\times \mathbb{C}^m$ making the following diagram commutative:
	\begin{equation}\label{diagram 3}
	\begin{tikzpicture}[descr/.style={fill=white,inner 		
		sep=2.5pt}]
	\matrix (m) [matrix of math nodes, row sep=3em,
	column sep=3em]
	{  Z^{\log}& &Q^{\log}_0\times \mathbb{C}^m\\
		&Q^{\log}_0& \\};
	\path[->,font=\scriptsize]
	(m-1-1) edge node[above] {$\iota^{\loga}$} (m-1-3) 
	(m-1-1) edge node[left] {$p^{\loga}$} (m-2-2) 
	(m-1-3) edge node[right] {$\pi_1$} (m-2-2);
	\end{tikzpicture}.
	\end{equation}
	By Proposition $\ref{prop:log triangularly smoothness}$, the canonical projection \[\pi_{Z^{\log}}:\mathcal{Q}_{Z^{\log}}(\mathfrak{X}_{Z^{\log}},\mathcal{M}_{\mathfrak{X}_{Z^{\log}}})\to Z^{\log}\] is smooth in a neighborhood of $(z_{q_0^{\loga}},q_0^{\loga},q_0^{\loga})$,  hence we can  opportunely choose (see \cite[p. 35]{Stie}, \cite[p. 269]{Lepot} and \cite[p. 94]{Raf}) a local trivialization 
	\begin{equation}\label{gamma log}
	\gamma^{\loga}:Z^{\log}\times Q^{\log}_0\to \mathcal{Q}_{Z^{\log}}(\mathfrak{X}_{Z^{\log}},\mathcal{M}_{\mathfrak{X}_{Z^{\log}}}).
	\end{equation}
	We notice that  the restriction of the canonical relative log cuirasse $\mathfrak{q}^{\loga}$ \eqref{log can} to $Z^{\log}$ produces a canonical relative log cuirasse on $(\mathfrak{X}_{Z^{\log}},\mathcal{M}_{\mathfrak{X}_{Z^{\log}}})$ over $Z^{\log}$. Hence,
	by $\eqref{def:map psi simple}$, we get a map
	\begin{equation}\label{eq:psi zlog simpler}
	\psi_{\mathfrak{q}^{\loga}}:\mathcal{Q}_{Z^{\log}}(\mathfrak{X}_{Z^{\log}},\mathcal{M}_{\X_{Z^{\log}}})\to Z^{\log}.
	\end{equation}
	Thus, we can define morphisms:
	\begin{equation}\label{def:morph omega log}
	\omega^{\loga}:=\psi_{\mathfrak{q}^{\loga}}\circ \gamma^{\loga} \text{ and } \delta^{\loga}:=\omega^{\loga}|_{\{(z_{q^{\loga}_0},q^{\loga}_0)\}\times Q^{\log}_0}
	\end{equation}
	By Proposition \ref{prop:log rel finite}, we can draw the following commutative diagram:
	\begin{equation}\label{diagram 4}
	\begin{tikzpicture}[descr/.style={fill=white,inner 		
		sep=2.5pt}]
	\matrix (m) [matrix of math nodes, row sep=3em,
	column sep=3em]
	{  Q^{\log}_0&Z^{\log}& &Q^{\log}_0\times \mathbb{C}^m\\
		& &Q^{\log}_0& \\};
	\path[->,font=\scriptsize]
	(m-1-2) edge node[above] {$\iota^{\loga}$} (m-1-4) 
	(m-1-2) edge node[left] {$p^{\loga}$} (m-2-3) 
	(m-1-1) edge node[above] {$\delta^{\loga}$} (m-1-2) 
	(m-1-4) edge node[right] {$\pi_1$} (m-2-3);
	\end{tikzpicture}.
	\end{equation}
	Proposition \ref{prop:log finite codim} is the log version of \cite[p. 595, Proposition 2]{DouVers} and \cite[p. 40, Satz 1.31]{Stie}.  See also \cite[p. 77, Proposition 4.40]{Raf}.
	\begin{proposition}\label{prop:log finite codim}
		The linear tangent map:
		\[\text{T}_{q^{\loga}_0}(p^{\loga}\circ\delta^{\loga}):\text{T}_{q^{\loga}_0}Q^{\log}_0\rightarrow \text{T}_{q^{\loga}_0}Q^{\log}_0\]
		is of the form $\Id-v^{\loga}$, with $v^{\loga}$ compact.
	\end{proposition}
	From Proposition $\ref{prop:log finite codim}$, it follows that $\ker \text{T}_{q^{\loga}_0}(p^{\loga}\circ \delta^{\loga})$ is of finite dimension. 
	Moreover, by \eqref{diagram 4}, we have
	\[\ker \text{T}_{q^{\loga}_0}(p^{\loga}\circ \delta^{\loga})\supset \ker \text{T}_{q^{\loga}_0}\delta^{\loga}=\ker \text{T}_{q^{\loga}_0}(\iota^{\loga}\circ\delta^{\loga}).\]
	Hence, $\ker \text{T}_{q^{\loga}_0}(\iota^{\loga}\circ\delta^{\loga})$ is of finite dimension.  Since $\pi_1$ \eqref{diagram 4} is a surjective map, we can conclude that $\Ima \text{T}_{q^{\loga}_0}(\iota^{\loga}\circ \delta^{\loga})$ has finite codimension in $\text{T}_{q^{\loga}_0}Q_0^{\log}$ (see \cite[p. 45]{Stie}). 
	
	Let us consder $\omega^{\loga}:Z^{\log}\times Q^{\log}_0\to Z^{\log}$ given by \eqref{def:morph omega log}.  
	Proposition $\ref{prop:log versality finite dim}$  is the log version of \cite[p. 36, Satz 1.25]{Stie}. The same proof applies here likewise. For further details see \cite[p. 94, Proposition 5.28]{Raf}.  
	
	\begin{proposition}\label{prop:log versality finite dim}
		Let $S$ be a Banach analytic space and $f,g:S\to Z^{\log}$ morphisms. Then $f^{*}(\X_{Z^{\log}},\mathcal{M}_{\X_{Z^{\log}}})\simeq g^{*}(\X_{Z^{\log}},\mathcal{M}_{\X_{Z^{\log}}})$, if and only if there exists $h^{\loga}:S\to Q^{\log}_0$ such that the following diagram commutes
		\begin{equation}\label{diagram log}
		\begin{tikzpicture}[descr/.style={fill=white,inner 		
			sep=2.5pt}]
		\matrix (m) [matrix of math nodes, row sep=3em,
		column sep=3em]
		{ & Z^{\log}\times Q^{\log}_0\\
			S&\\
			& Z^{\log} \\};
		\path[->,font=\scriptsize]
		(m-2-1) edge node[left] {$(f,h^{\loga})$} (m-1-2) 
		(m-2-1) edge node[left] {$g$} (m-3-2) 
		(m-1-2) edge node[right] {$\omega^{\loga}$} (m-3-2);
		\end{tikzpicture}.
		\end{equation}
	\end{proposition}
	In other words, \[f^{*}(\X_{Z^{\log}},\mathcal{M}_{\X_{Z^{\log}}})\simeq g^{*}(\X_{Z^{\log}},\mathcal{M}_{\X_{Z^{\log}}})\] if and only if, for each $s\in S$, $g(s)$ is obtained \enquote{changing} $f(s)$ by a log cuirasse $q^{\loga}$ on the central fibre $(X_0,\mathcal{M}_{X_0})$. Notice that, by Proposition \ref{prop:log versality finite dim}
	\begin{equation}
	\omega^{\loga}|_{Z^{\log}\times q_0^{\loga}}=\Id_{Z^{\log}}.
	\end{equation}

	Let us denote with $\Ex^1(X_0,\mathcal{M}_{X_0})$ the set of equivalence classes of infinitesimal deformations of $(X_0,\mathcal{M}_{X_0})$, that is deformations over the double point $D=(\{\cdot\},\mathbb{C}[\epsilon]/\epsilon^{2})$. For the sake of clarity, set $r^{\loga}_0:=(z_{q^{\loga}_{0}},q^{\loga}_{0})$.
	
	Since $(\mathfrak{X}_{Z^{\log}},\mathcal{M}_{\mathfrak{X}_{Z^{\log}}})\to (Z^{\log},r^{\loga}_0)$ is complete, the Kodaira-Spencer map (\ref{ks}) is surjective
	\[\text{ks}:\text{T}_{r^{\loga}_0}Z^{\log}\twoheadrightarrow \Ex^{1}(X_0,\mathcal{M}_{X_0}).\]
	The kernel $\ker \text{ks}$ corresponds to the trivial deformations of $(X_0,\mathcal{M}_{X_0})$ over $D$. 
	By Proposition $\ref{prop:log versality finite dim}$, with $S=D$, we see that the trivial deformations of $(X_0,\mathcal{M}_{X_0})$ over $D$ are given by $\Ima \text{T}_{q^{\loga}_0}\delta^{\loga}$. Hence, 
	\begin{equation}
	\Ex^1(X_0,\mathcal{M}_{X_0})=\text{T}_{r^{\loga}_0}Z^{\log}/\Ima \text{T}_{q^{\loga}_0}\delta^{\loga}.
	\end{equation}
	Let us identify $Z^{\log}$ with its image  in $Q^{\log}_0\times \mathbb{C}^{m}$ under $\iota^{\loga}$. By Proposition \ref{prop:log finite codim}, let $\Sigma^{\log}$ be the Banach submanifold of $Q_0^{\log}$ such that 
	\begin{equation}\label{sigmalog}
	\text{T}_{q^{\loga}_0}\Sigma^{\log}\oplus\ker\text{T}_{q^{\loga}_0}\delta^{\loga}=\text{T}_{q^{\loga}_0}Q^{\log}_0.
	\end{equation}
	Let $r:Q^{\log}_0\times \mathbb{C}^m\to \delta^{\loga}(\Sigma^{\log})$ be a retraction and set
	\begin{equation}\label{Rlog}
	R^{\log}:=r^{-1}(q^{\loga}_0)\cap Z^{\log}.
	\end{equation}
	By construction
	\begin{equation}\label{effect}
	\text{T}_{r^{\loga}_0}R^{\log}=\Ex^1(X_0,\mathcal{M}_{X_0})
	\end{equation}
	
	\begin{lemma}(\cite[p. 598, Proposition 5]{DouVers} and \cite[p. 37, Satz 1.28]{Stie})\label{sandwich}
		Let $\Sigma_1, H$ and $\Sigma_2$ be Banach manifolds, with $\Sigma_1, H\subset \Sigma_2$ and  $0\in \Sigma_1\cap H$. Assume that $\Sigma_1$ is of finite codimension and
		\begin{equation}
		\text{T}_0\Sigma_1\oplus \text{T}_0 H =\text{T}_0\Sigma_2.
		\end{equation}
		Let $Y$ be another subspace of $\Sigma_2$, containing $\Sigma_1$, and set
		\[R:=H\cap Y.\]
		Let 
		\[\phi:\Sigma_1\times R\to Y\]
		be a morphism inducing the identity on $\Sigma_1\times 0$ and $0\times R$. Then, $\phi$ is an isomorphism.
	\end{lemma}
	From Lemma \ref{sandwich}, we obtain that the restriction of the morphism \eqref{def:morph omega log}
	\begin{equation}\label{def:morph omega l}
	\omega^{\loga}|_{R^{\log}\times\Sigma^{\log}}:R^{\log}\times\Sigma^{\log}\to Z^{\log} 
	\end{equation} 
	is an isomorphism. This fact, together with Proposition \ref{prop:log versality finite dim} and \eqref{effect}, is used to prove Theorem \ref{thm:log fin dim douady versal}.
	\smallskip
	
	Let $i:R^{\log}\hookrightarrow Z^{\log}$ be the canonical injection. Set
	\[(\mathfrak{X}_{R^{\log}},\mathcal{M}_{\mathfrak{X}_{R^{\log}}}):=i^{*}(\mathfrak{X}_{Z^{\log}},\mathcal{M}_{\mathfrak{X}_{Z^{\log}}}).\]
	Let $\alpha_{q^{\loga}_0}:(\mathfrak{X}_{R^{\log}},\mathcal{M}_{\mathfrak{X}_{R^{\log}}})(r^{\loga}_0)\rightarrow (X_0,\mathcal{M}_{X_0})$ be the log isomorphism given by Proposition $\ref{prop:log iso alpha}$. 	
	The proof of Theorem \ref{thm:log fin dim douady versal} is identical to the proof of \cite[p. 598, Théorème Principal and p. 601, Proposition 1]{DouVers} and  to the proof of  \cite[p. 38, Satz 1.30]{Stie}. 
	\begin{theorem}\label{thm:log fin dim douady versal}
		The triple $((R^{\log},r^{\loga}_0), (\mathfrak{X}_{R^{\log}},\mathcal{M}_{\mathfrak{X}_{R^{\log}}}), \alpha_{q^{\loga}_0})$, is a semi-universal  deformation of $(X_0,\mathcal{M}_{X_0})$.
	\end{theorem}
	
	\proof
	Let $((S,s_0),(X,\mathcal{M}_X),i)$ be a deformation of $(X_0,\mathcal{M}_{X_0})$. By Proposition $\ref{prop:versal Z}$,  $((\mathfrak{X}_{Z^{\log}},\mathcal{M}_{\mathfrak{X}_{Z^{\log}}})\to Z^{\log},(z_{q^{\loga}_0},q^{\loga}_0))$ is a complete deformation of of $(X_0,\mathcal{M}_{X_0})$. Hence, there exists a morphism $\psi_{\loga}:S\to Z^{\log}$ such that \[(X,\mathcal{M}_X)\simeq \psi_{\loga}^{*}(\mathfrak{X}_{Z^{\log}},\mathcal{M}_{\mathfrak{X}_{Z^{\log}}}).\] 
	Let $\Sigma^{log}$ and $R^{\log}$ given by \eqref{sigmalog} and \eqref{Rlog} respectively. Let $\pi_{R^{\log}}:R^{\log}\times\Sigma^{\log}\to R^{\log}$ and $\pi_{\Sigma^{\log}}:R^{\log}\times\Sigma^{\log}\to \Sigma^{\log}$ be the projections. By Lemma $\ref{sandwich}$, the morphism $\omega^{\loga}|_{R^{\log}\times\Sigma^{\log}}:R^{\log}\times\Sigma^{\log}\to Z^{\log}$  $\eqref{def:morph omega l}$ is an isomorphism. Thus, setting $g:=\pi_{R^{\log}}\circ (\omega^{\loga}|_{R^{\log}\times\Sigma^{\log}})^{-1}\circ \psi_{\loga}$ and $h^{\loga}:=\pi_{\Sigma^{\log}}\circ (\omega^{\loga}|_{R^{\log}\times\Sigma^{\log}})^{-1}\circ \psi_{\loga}$, we have
	\[\omega^{\loga}\circ(g,h^{\loga})=\psi_{\loga}.\]
	Hence, by Proposition $\ref{prop:log versality finite dim}$
	\[g^{*}(\mathfrak{X}_{R^{\log}},\mathcal{M}_{\mathfrak{X}_{R^{\log}}})\simeq \psi_{\loga}^{*}(\mathfrak{X}_{Z^{\log}},\mathcal{M}_{\mathfrak{X}_{Z^{\log}}})\simeq (X,\mathcal{M}_X).\]
	Moreover, by construction, $T_{q_0^{\loga}}R^{\log}=\Ex^1(X_0,\mathcal{M}_{X_0})$ \eqref{effect}.  
	
	Thus, the deformation $((\mathfrak{X}_{R^{\log}},\mathcal{M}_{\mathfrak{X}_{R^{\log}}})\rightarrow R^{\log},r^{\loga}_0)$ is complete and effective.
	\smallskip 
	
  	Now, let $((S,s_0),(X,\mathcal{M}_X),i)$ be a deformation of $(X_0,\mathcal{M}_{X_0})$ and $(S',s_0)$ a subgerm of $(S,s_0)$. Because of the just proved completeness, we can find a morphism $h':(S',s_0)\to (R^{\log},r_0)$ such that \[(X,\mathcal{M}_X)|_{S'}\simeq h^{'*}(\mathfrak{X}_{R^{\log}},\mathcal{M}_{\mathfrak{X}_{R^{\log}}}).\] Let $\mathfrak{q}^{\dagger}$ be the canonical relative log cuirasse on $(\mathfrak{X}_{Z^{\log}},\mathcal{M}_{\mathfrak{X}_{Z^{\log}}})$ over $Z^{\log}$ $\eqref{log can}$. Then, $h^{'*}\mathfrak{q}^{\dagger}$ is a relative log cuirasse on $(X,\mathcal{M}_{X})|_{S'}$ over $S'$, whose associated morphism $\eqref{log can cuirasse}$ coincides with $h'$. Since, by Proposition $\ref{prop:log triangularly smoothness}$, $\mathcal{Q}_S(X,\mathcal{M}_X)$ is smooth over $S$ in a neighborhood of $q^{\dagger}_0\in\mathcal{Q}((X(s_0),\mathcal{M}_{X(s_0)})$,  there exists a relative cuirasse $q^{\dagger}$ on $(X,\mathcal{M}_X)$ over $S$, such that $q^{\dagger}|_{S'}=h^{'*}\mathfrak{q}^{\dagger}$. Let $\tilde{h}:S\to Z^{\log}$ be the morphism associated to $q^{\dagger}$ $\eqref{def:map psi simple}$ and $\pi_{R^{\log}}:Z^{\log}\to R^{\log}$ the projection. Then, $h:=\pi_{R^{\log}}\circ \tilde{h}$ satisfies $(X,\mathcal{M}_X)\simeq h^{*}(\mathfrak{X}_{R^{\log}},\mathcal{M}_{\mathfrak{X}_{R^{\log}}})$ and $h|_{S'}=h'$. 
	
	Thus, the deformation $(\mathfrak{X}_{R^{\log}},\mathcal{M}_{\mathfrak{X}_{R^{\log}}})\rightarrow R^{\log}$ is also versal and, therefore,	semi--universal.
	\endproof
	
	\section{Semi-universal deformations of log morphisms}
	\label{log morph}

	In what follows, we  construct a semi-universal deformation of a  morphism $f_0:(X_0,\mathcal{M}_{X_0})\rightarrow (Y_0,\mathcal{M}_{Y_0})$ of compact fine log complex spaces.
	Let $X$ be a complex space and $\alpha_i:\M_i\to\O_X$, $i=1,2$,  two fine log structures on $X$. Let $\gamma:\overline{\M}_1\to\overline{\M}_2$ be a morphism of the ghost sheaves. Let $f:T\to X$ be a morphism of complex spaces and set $\gamma_T:(\overline{\M}_1)_T\to(\overline{\M}_2)_T$, the pull-back of $\gamma$ via $f$.

	\begin{lemma}\label{lem:siebert}(\cite[p. 474]{GSGW})
		The functor
		\[\Mor^{\log}_{X}: \textit{An}_{X}\to \textit{Sets}\]
		defined on the objects by
		\[(f:T\to X)\mapsto\{\varphi:(T,f^{*}\M_1)\to(T,f^{*}\M_2)|\overline{\varphi}^{\flat}=\gamma_T\}\]
		is represented by a complex space $\mathfrak{Mor}^{\log}_X(\mathcal{M}_1,\mathcal{M}_2)$ over $X$.
	\end{lemma}
	\proof
	By the universal property, the statement is local in $X$. Hence, let $\beta_i:P_i\to\Gamma(X,\M_i)$, $i=1,2$, be two log charts for $\M_1$ and $\M_2$ respectively. Let $p_1,...,p_n\in P_1$ be a generating set for $P_1$ as monoid. Consider the sheaf of finitely generated $\O_X$-algebras
	\[\mathcal{F}_X:=\O_X[P^{gp}_1]/\langle\alpha_1(\beta_1(p_i))-z^{p_i}\alpha_2(\beta_2(\gamma(p_i)))|1\leq i\leq n\rangle.\]
	Set $\mathfrak{Mor}^{\log}_X(\mathcal{M}_1,\mathcal{M}_2):=\boldsymbol{\Specan } \mathcal{F}_X$, the relative analytic spectrum of $\mathcal{F}_X$ over $X$. Now, we check the universal property. Let $f:T\to X$ be given. We want to show that giving a commutative diagram of complex spaces
	
	\[
	\begin{tikzpicture}[descr/.style={fill=white,inner 		
		sep=2.5pt}]
	\matrix (m) [matrix of math nodes, row sep=3em,
	column sep=3em]
	{T &    \boldsymbol{\Specan } \mathcal{F}_X\\
		& X \\ };
	\path[->,font=\scriptsize]
	(m-1-1) edge node[above] {$g$} (m-1-2) 
	(m-1-1) edge node[below] {$f$} (m-2-2)
	(m-1-2) edge node[right] {$ $} (m-2-2);
	\end{tikzpicture}
	\]
	is equivalent to giving a log morphism $\varphi:(T,f^{*}\M_1)\to (T,f^{*}\M_2)$, which is the identity on $X$ and such that $\overline{\varphi}^{\flat}=\gamma_T$.
	Giving a morphism $g$ is equivalent to giving a section of $(\boldsymbol{\Specan}\mathcal{F}_X)\times_X T$ over $T$. But
	\[(\boldsymbol{\Specan}\mathcal{F}_X)\times_X T=\boldsymbol{\Specan}\O_T[P^{gp}_1]/\langle f^{*}(\alpha_1(\beta_1(p_i)))-z^{p_i}f^{*}(\alpha_2(\beta_2(\gamma(p_i))))|1\leq i\leq n\rangle,\]
	and the latter complex space is $\boldsymbol{\Specan }\mathcal{F}_T$ associated to the data $(T,f^{*}\M_1)$,$(T,f^{*}\M_2)$ with charts $f^{*}(\beta_i)=f^{\flat}\circ\beta_i:P_i\to\Gamma(Y,f^{*}\M_i)$. Thus, without loss of generality, we can assume $T=X$ and $f$ is the identity. Now, giving $\varphi:(X,\M_1)\to (X,\M_2)$, with $\overline{\varphi}^{\flat}=\gamma$, is equivalent to specifying $\varphi^{\flat}$. From $\varphi^{\flat}$ we obtain a map $\eta:P_1\to\Gamma(X,\O^{\times}_X)$ with the property that for all $p\in P_1$,
	\[\varphi^{\flat}(\beta_1(p))=\eta(p)\cdot\beta_2(\gamma(p)).\]
	Conversely, $\eta$ completely determines $\varphi^{\flat}$. In addition, $\varphi^{\flat}$ is a homomorphism of monoids if and only if $\eta$ is a homomorphism, and since $\eta$ takes values in the group $\O^{\times}_X$, specifying $\varphi^{\flat}$ is equivalent to specifying a section of $\boldsymbol{\Specan}\mathcal{O}_X[P^{gp}_1]$. Indeed, a section of $\boldsymbol{\Specan }\mathcal{O}_X[P^{gp}_1]$ over $X$ is the same as a morphism $X\to \Spec \mathbb{C}[P^{gp}_1]$, which in turn is the same as an element of $\Hom(P_1,\Gamma(X,\O^{\times}_X))$. Second, since $\varphi^{*}=\id$, we must have $\alpha_1=\alpha_2\circ\varphi^{\flat}$, so for each $p\in P_1$, we must have
	\[\alpha_1(\beta_1(p))=\alpha_2(\varphi^{\flat}(\beta_1(p)))=\eta(p)\cdot \alpha_2(\beta_2(\gamma(p))).\]
	If this holds for each $p_i$, it holds for all $p$. Thus a section of $\boldsymbol{\Specan }\O_X[P^{gp}_1]$ over $X$ determines a morphism of log structures if and only if it lies in the subspace determined by the equations 
	\[\alpha_1(\beta_1(p_i))-z^{p_i}\alpha_2(\beta_2(\gamma(p_i))),\]
	demonstrating the result.
	\endproof
	
	Now, assume the complex space $X$ is proper over a germ of complex spaces $(S,s_0)$.
	
	\begin{proposition}\label{prop:siebert}(\cite[p. 475]{GSGW}) The functor \[\Mor^{\log}_{X/S}:(f:(T,t_0)\to (S,s_0))\mapsto\{\varphi:(X_T,(\M_1)_T)\to(X_T,(\M_2)_T)|\overline{\varphi}^{\flat}=\gamma_T\}\]
		is represented by a germ $\mathfrak{Mor}^{\log}_{X/S}(\mathcal{M}_1,\mathcal{M}_2)$ of complex spaces over $(S,s_0)$.
	\end{proposition}
	\proof
	Let $Z=\mathfrak{Mor}^{\log}_X(\mathcal{M}_1,\mathcal{M}_2)$. By Lemma $\ref{lem:siebert}$, $\mathfrak{Mor}^{\log}_{X/S}(\mathcal{M}_1,\mathcal{M}_2)$ is isomorphic to the functor
	\[(\psi:T\to S)\mapsto\{\text{sections of }\psi^{*}Z\to\psi^{*}X\}.\]
	This is exactly the functor of sections $\prod_{X/S}(Z/X)$ discussed, in the algebraic-geometric setting, in \cite[p. 267]{Gt} and here it is represented by an open subspace of the relative Douady space of $Z$ over $S$ (see \cite{PourRel}).
	
	\endproof
	
	\begin{proposition}\label{thm:universal complex space}(\cite[p. 130]{Flenner})
		Let $f_0:X_0\rightarrow Y_0$ be a holomorphic map between compact complex spaces. Then $f_0$ admits a semi-universal deformation.
	\end{proposition}
	
	\begin{theorem}\label{thm:vers log morph}
		Every morphism $f_{0}:(X_0,\M_{X_0})\to (Y_0,\M_{Y_0})$ of compact fine log complex spaces admits a semi-universal deformation parametrized by a germ of complex spaces $(S,s_0)$.
	\end{theorem}	
	\proof
	Let $((\mathfrak{X},\mathcal{M}_{\mathfrak{X}})\rightarrow R,r_0)$ and $((\mathcal{Y},\mathcal{M}_{\mathcal{Y}})\rightarrow R,r_0)$ be the semi-universal deformations of $(X_0,\mathcal{M}_{X_0})$ and $(Y_0,\mathcal{M}_{Y_0})$ respectively given by Theorem $\ref{thm:log fin dim douady versal}$. By pulling-back to the product of the base spaces, we can assume that the two deformations are defined over the same base space. Let us consider the finite dimensional complex analytic space $\mathfrak{Mor}_{R}(\mathfrak{X},\mathcal{Y})$
	given by Proposition $\ref{thm:universal complex space}$. Let $p:\mathfrak{Mor}_{R}(\mathfrak{X},\mathcal{Y})\to R$ be the projection and set $\underline{m}_0:=(r_0,\underline{f}_0)$. By Proposition $\ref{thm:universal complex space}$, we get a universal morphism $\underline{f}:p^{*}\X\to p^{*}\mathcal{Y}$, such that the restriction of $\underline{f}$ to the central fibre $p^{*}\X(\underline{m}_0)$ equals $\underline{f}_0$. 
	We can consider two fine log structures on $p^{*}\X$, namely \[\M_1:=p^{*}\M_{\X} \text{ and } \M_2:=\underline{f}^{*}p^{*}\mathcal{M}_{\mathcal{Y}}.\] Set
	$\gamma:=\overline{f}^{\flat}_0$ and  $m_0:=(\underline{m}_0,f^{\flat}_0)$. For the sake of clarity, denote $\mathfrak{M}:= \mathfrak{Mor}_{R}(\mathfrak{X},\mathcal{Y})$. 
	Now, consider the germ of complex spaces
	$(\mathfrak{Mor}^{\log}_{p^{*}\X/\mathfrak{M}}(\mathcal{M}_1,\mathcal{M}_2),m_0)$, together with the projection \[\pi:(\mathfrak{Mor}^{\log}_{p^{*}\X/\mathfrak{M}}(\mathcal{M}_1,\mathcal{M}_2),m_0)\to (\mathfrak{Mor}_{R}(\mathfrak{X},\mathcal{Y}),\underline{m}_0),\] provided by Proposition $\ref{prop:siebert}$. Moreover, by Proposition $\ref{prop:siebert}$, we get a morphism $f^{\flat}:\pi^{*}(p^{*}\X,\M_1)\to\pi^{*}(p^{*}\X,\M_2)$. 
	
	Hence, we get a log morphism \[f:=(\pi^{*}\underline{f},f^{\flat}):\pi^{*}p^{*}(\X,\M_{\X})\to\pi^{*}p^{*}(\mathcal{Y},\M_{\mathcal{Y}})\] over $(\mathfrak{Mor}^{\log}_{p^{*}\X/\mathfrak{M}}(\mathcal{M}_1,\mathcal{M}_2),m_0)$. Set $(S,s_0):=(\mathfrak{Mor}^{\log}_{p^{*}\X/\mathfrak{M}}(\mathcal{M}_1,\mathcal{M}_2),m_0)$. Using the universal property of $(S,s_0)$ (see Proposition $\ref{prop:siebert}$) and Theorem $\ref{thm:log fin dim douady versal}$, the statement follows. 
	\endproof
	Moreover, we can deform $(X_0,\mathcal{M}_{X_0})$ as relative log space over $(Y_0,\mathcal{M}_{Y_0})$. That is, we can deform $(X_0,\mathcal{M}_{X_0})$ together with the morphism $f_0$ into $(Y_0,\mathcal{M}_{Y_0})$. In this case, $Y_0$ needs not to be compact. More precisely,
	\begin{definition}\label{def:relative 1}
		Let $f_0:(X_0,\mathcal{M}_{X_0})\to(Y_0,\mathcal{M}_{Y_0})$ be a log  morphism of fine log complex spaces, with $X_0$ compact.
		A semi-universal deformation of $(X_0,\mathcal{M}_{X_0})$ \textit{over} $(Y_0,\mathcal{M}_{Y_0})$, with base a germ of complex spaces $(S,s_0)$, is a commutative diagram 
		\[	\begin{tikzpicture}[descr/.style={fill=white,inner 		
			sep=2.5pt}]
		\matrix (m) [matrix of math nodes, row sep=3em,
		column sep=3em]
		{ (\X,\mathcal{M}_{\X}) &   & (Y_0,\mathcal{M}_{Y_0})\times S\\
			& (S,s_0)  & \\ };
		\path[->,font=\scriptsize]
		(m-1-1) edge node[above] {$f$} (m-1-3) 
		(m-1-1) edge node[left] {$p$} (m-2-2)
		(m-1-3) edge node[right] {$\pi_2$} (m-2-2);
		\end{tikzpicture},\] 
		where $p$ is a semi-universal deformation of $(X_0,\mathcal{M}_{X_0})$, together with an isomorphism $i:(X_0,\mathcal{M}_{X_0})\to (\X,\mathcal{M}_{\X})(s_0)$, such that $f\circ i=f_0$.
	\end{definition}
	
	The same proof of Theorem $\ref{thm:vers log morph}$, with $(\mathcal{Y},\mathcal{M}_{\mathcal{Y}}):=(Y_0,\mathcal{M}_{Y_0})\times R$, gives us the following 
	
	\begin{corollary}\label{cor:def log morph}
		Let $f_0:(X_0,\mathcal{M}_{X_0})\to(Y_0,\mathcal{M}_{Y_0})$ be a log  morphism of fine log complex spaces, with $X_0$ compact. Then $(X_0,\mathcal{M}_{X_0})$ admits a semi-universal deformation over $(Y_0,\mathcal{M}_{Y_0})$.
	\end{corollary}
	
	\begin{remark}\label{rem:relative 3}
		If $f_0$ is a log embedding, then Corollary $\ref{cor:def log morph}$ gives us a semi-universal deformation of a log subspace $(X_0,\mathcal{M}_{X_0})$ in a fixed ambient log space $(Y_0,\mathcal{M}_{Y_0})$.
	\end{remark}

	Now, we assume that $f_0$ is a log flat (log smooth) morphism. We show that, in this case, we get a log flat (log smooth) semi-universal deformation of $f_0$.	 We need the following two results in analytic geometry.

	\begin{proposition}(\textit{Critère de platitude par fibres})\label{Gro}(\cite[p. 138]{EGA43} and \cite[\href{https://stacks.math.columbia.edu/tag/00MP}{Tag 00MP}]{stacks-project})
		Let $S$ be a complex space.
		Let $f : X \to Y$ be a morphism of complex spaces over $S$.
		Let $\mathcal{F}$ be a coherent $\mathcal{O}_X$-module.
		Let $x \in X$. Set $y = f(x)$ and $s \in S$ the image of $x$ in $S$. Let $X_s$ and $Y_s$ be the fibres of $X$ and $Y$ over $s$. Set:
		\[\mathcal{F}_s = (X_s \hookrightarrow X)^*\mathcal{F}.\]
		Assume $\mathcal{F}_x \not = 0$.
		Then the following are equivalent:
		\begin{enumerate}
			\item $\mathcal{F}$ is flat over $S$ at $x$ and
			$\mathcal{F}_s$ is flat over $Y_s$ at $x$;
			\item $Y$ is flat over $S$ at $y$ and $\mathcal{F}$ is
			flat over $Y$ at $x$.
		\end{enumerate}
	\end{proposition}

	\begin{proposition}\label{log smoothness}(\cite[p. 159]{Fischer1976})
		Let $f:X\rightarrow Y$ be a morphism of complex spaces. Let $p\in X$. Then the following are equivalent
		\begin{enumerate}
			\item $f$ is smooth (submersion) at $p\in X$;
			\item $f$ is flat at $p$ and the fibre $X_{f(p)}$ is a manifold.
		\end{enumerate}
	\end{proposition}
Proposition \ref{Gro} is due to A. Grothendieck in the algebraic geometry setting.  The result can be naturally extended to the analytic setting as for any complex space $(X,\mathcal{O}_X)$ and  $p\in X$, the stalk $\mathcal{O}_{X,p}$ is a Noetherian local ring (see \cite[p. 80]{kaup}). 
	
	\begin{lemma}\label{closed maps}
		Let $f: X\rightarrow Y$ be a continuous map between topological spaces. If $f$ is
		closed, then for all $y\in Y$ and open subset $U\subset X$ satisfying $f^{-1}(y)\subset U$, there exists an
		open neighborhood $V$ of $y$ satisfying $f^{-1}(V)\subset U$.
	\end{lemma}
	\proof
	Let us consider the closed subset $X\backslash U$. Since $f$ is closed, $f(X\backslash U)$ is closed in $Y$. Therefore, $Y\backslash f(X\backslash U)$ is open in $Y$ and it contains $y$ as $f^{-1}(y)\subset U$. Take $V:=f^{-1}(Y\backslash f(X\backslash U)).$ 
	\endproof
	
The following Lemma \ref{prop:log smooth log flat} can be found, in the algebraic geometry setting, in \cite[p. 424]{Ogus18}. This is a local statement, which extends naturally to the analytic setting. 
	\begin{lemma}\label{prop:log smooth log flat}
		Any log smooth morphism of fine log complex spaces is log flat.
	\end{lemma}
	
	Let $f:(\X,\mathcal{M}_{\X}) \to(\mathcal{Y},\mathcal{M}_{\mathcal{Y}})$ be the semi-universal deformation of $f_{0}:(X_0,\M_{X_0})\to (Y_0,\M_{Y_0})$, over a germ of complex spaces $(S,s_0)$, given by Theorem $\ref{thm:vers log morph}$ or Corollary $\ref{cor:def log morph}$. Denote with $\pi_1$ and $\pi_2$ the morphisms of $(\X,\mathcal{M}_{\X})$ and $(\mathcal{Y},\mathcal{M}_{\mathcal{Y}})$ into $(S,s_0)$ respectively.
	\begin{proposition}\label{prop:log flatness}
		If $f_0$ is log flat (log smooth), then  $f$ is log flat (log smooth) in an open neighborhood of $s_0$.
	\end{proposition}
	\proof
	Let us assume that there exists an open neighborhood $U'$ of $X_0$ in $\X$ such that $f|_{(U',\mathcal{M}_{U'})}$ is log flat (log smooth). Then, since $\pi_1:\X\to S$ is a proper map between locally compact Hausdorff spaces, it is closed. Hence, by Lemma $\ref{closed maps}$, we can find an open neighborhood $W$ of $s_0$ such that  $\pi_1^{-1}(W)$ is contained in $U'$. This ensures us that $f$ is log flat (log smooth) as relative morphism over $(W,s_0)\subset (S,s_0)$. Since log flatness (log smoothness) is a local property, we choose a log chart for $f$. We have the following commutative diagram
	\[	\begin{tikzpicture}[descr/.style={fill=white,inner 		
		sep=2.5pt}]
	\matrix (m) [matrix of math nodes, row sep=3em,
	column sep=3em]
	{ U &  \Spec\mathbb{C}[P]\\
		V & \Spec\mathbb{C}[Q] \\ };
	\path[->,font=\scriptsize]
	(m-1-1) edge node[above] {$\beta$} (m-1-2)
	(m-1-1) edge node[left] {$ \underline{f} $} (m-2-1)
	(m-1-2) edge node[right] {$g$} (m-2-2)
	(m-2-1) edge node[above] {$\gamma$} (m-2-2);
	\end{tikzpicture}.\]
	Let us consider the universal morphism $\underline{u}:U\rightarrow V\times_{\Spec\mathbb{C}[Q]} \Spec\mathbb{C}[P]$. Let $p:V\times_{\Spec\mathbb{C}[Q]} \Spec\mathbb{C}[P]\rightarrow V$ be the projection. 
	 We get the following commutative diagram:
	\[	\begin{tikzpicture}[descr/.style={fill=white,inner 		
		sep=2.5pt}]
	\matrix (m) [matrix of math nodes, row sep=3em,
	column sep=3em]
	{ U &   & V\times_{\Spec\mathbb{C}[Q]} \Spec\mathbb{C}[P]\\
		& (S,s_0)  & \\ };
	\path[->,font=\scriptsize]
	(m-1-1) edge node[above] {$\underline{u}$} (m-1-3) 
	(m-1-1) edge node[left] {$\pi_1|_{U}$} (m-2-2)
	(m-1-3) edge node[right] {$\pi_2\circ p$} (m-2-2);
	\end{tikzpicture}.\] 
	Assume $f_0$ log flat, then  we have that $\underline{u}$ is flat at $s_0$. Moreover, by Theorem $\ref{thm:log fin dim douady versal}$, $\pi_1|_{U}$ is flat too.  For the sake of clarity, set $A:=V\times_{\Spec\mathbb{C}[Q]} \Spec\mathbb{C}[P]$. 
	
	We use Proposition $\ref{Gro}$ for  $\mathcal{F}=\mathcal{O}_{U}$. Since \textit{condition} $1$ holds, by \textit{condition} $2$ we get that $\mathcal{O}_{U,x}$ is a flat $\mathcal{O}_{A,\underline{u}(x)}$-module, for each $x\in \pi_1^{-1}|_{U}(s_0)$. Since every flat holomorphic map is open, we get the existence of an open subset $U'$ of $U$, containing $\pi_1^{-1}|_{U}(s_0)$, such that $\underline{u}_{\arrowvert U'}$ is flat.  This proves the first part of the statement. Now, assume $f_0$ log smooth. By Lemma $\ref{prop:log smooth log flat}$, $f_0$ is log flat. Hence, by the first part of this proof, we get the existence of an open subset $U'$ in $U$ such that $\underline{u}_{\arrowvert U'}$ is flat. Let $x\in \pi_1^{-1}|_{U'}(s_0)$ and set $y:=\underline{u}(x)\in V\times_{\Spec\mathbb{C}[Q]} \Spec\mathbb{C}[P]$. Since $f_0$ is log smooth, $\underline{u}$ is smooth at $s_0$. Hence,  we get that the fibre $U'_{y}$ of $\underline{u}_{\arrowvert U'}$ over $y$ is a manifold. Therefore, using Proposition $\ref{log smoothness}$, we get the second part of the statement.
	
	\endproof

	\begin{example}
		Let $(\Spec\mathbb{C},Q)$ be a log point. Let $(X_0,\mathcal{M}_{X_0})$ be a compact fine log complex space and $f_0:(X_0,\mathcal{M}_{X_0})\rightarrow (\Spec\mathbb{C},Q)$ a log smooth morphism.  Since $\Hom((Q,+),(\mathbb{C},\cdot))=\Spec\mathbb{C}[Q]$, a semi-universal deformation of the log point $(\Spec\mathbb{C},Q)$ is given by the affine toric variety $\Spec\mathbb{C}[Q]$ endowed with the canonical log structure. Let $p_0\in \Spec\mathbb{C}[Q]$ be the base point. Let $((\X,\mathcal{M}_{\X})\to R,r_0)$ be the semi-universal deformation of $(X_0,\mathcal{M}_{X_0})$ given by Theorem $\ref{thm:log fin dim douady versal}$. 
		Let  $R\times \Spec\mathbb{C}[Q]$ and consider the projections $\pi_1$, $\pi_2$  onto the first and  second factor respectively. Then,  $\pi_1^{*}(\X,\mathcal{M}_{\X})$ and $\pi_2^{*}\Spec(Q\to\mathbb{C}[Q])$ are semi-universal deformations of $(X_0,\mathcal{M}_{X_0})$ and  $(\Spec\mathbb{C},Q)$  over $R\times \Spec\mathbb{C}[Q]$ respectively . Let $(r_0,p_0)\in R\times \Spec\mathbb{C}[Q]$ be the base point.
		By Theorem $\ref{thm:vers log morph}$, we get a germ of complex spaces $(S,s_0)$, together with a morphism of germs $p:(S,s_0)\to (R\times \Spec\mathbb{C}[Q], (r_0,p_0))$, and a log $S$-morphism $f:p^{*}\pi_1^{*}(\X,\mathcal{M}_{\X})\to p^{*}\pi_2^{*}\Spec(Q\to\mathbb{C}[Q])$, which is a semi-universal deformation of $f_0$. By Proposition $\ref{prop:log flatness}$, $f$ is log smooth in a neighborhood of $(X_0,\mathcal{M}_{X_0})$.
	\end{example}
\appendix
\section{Existence of directed log charts}
\label{appendix directed}
Let $(X,\M_{X})$ be a compact fine log complex space. We recall from Definition \ref{def:compatible log charts} that a set of \textit{directed} log charts is a collection of  log charts $(\theta_i:P_i\to\mathcal{M}_{U_i})_{i\in J}$ covering $(X,\mathcal{M}_{X})$, together with a morphism
\[\varphi^{i}_j:P_i\to P_j\oplus \mathcal{O}^{\times}_{U_j},\]
for each $j\in J_1\cup J_2$ and $i\in\partial j$, such that
\[(\theta_j\cdot \Id_{\mathcal{O}^{\times}_{U_j}})\circ \varphi^i_j=\theta_i|_{U_j}.\]
In what follows, we show that there exists a finite collection of directed log charts on $(X,\M_{X})$.

\begin{definition}\label{poset log charts}
	Let $\theta_i:P_i\to \Gamma(U_i,\M_{X})$, $i=1,2$, be two log charts.
	We write $\theta_2\leq\theta_1$ if and only if $U_2\subseteq U_1$ and there exists a morphism 
	\[\varphi:P_1\to P_2\oplus \mathcal{O}^{\times}_{U_2},\]
	such that
	\[(\theta_2\cdot \Id_{\mathcal{O}^{\times}_{U_2}})\circ \varphi=\theta_1|_{U_2}.\]
	Moreover, we say that $\theta_1$ and $\theta_2$ are \textit{equivalent} if and only if $\theta_1\leq\theta_2$ and $\theta_2\leq\theta_1$.
\end{definition}

Since $\M_X$ is a fine log structure, it induces a stratification of $X$ such that the restriction of $\ol \M_X$ to each stratum is locally constant. 

\begin{definition}\label{def:centered charts}
	Let 
	$\theta_i:P_i\to \Gamma(U_i,\M_X)$ be a chart and $x\in X$. 	We call $\theta_i$ \textit{a chart centered at} $x$, if it induces an isomorphism $P_i\simeq \Gamma(U_i,\ol\M_X)$ and  the restriction map
	$\Gamma(U,\ol \M_X)\to \ol \M_{X,x}$
	is an isomorphism.
	
\end{definition}

\begin{lemma}\label{lem:equivalent charts}
	Let $\theta_i:P_i\to\Gamma(U_i,\M_X)$, $i=1,2$, be two log charts with $U_2\subseteq U_1$. Let $Z_i\subset X$, $i=1,2$, be strata with $ Z_1\subseteq \ol Z_2$. Assume $\theta_1$ and $\theta_2$ are centered at points $x_1\in Z_1\cap U_1$ and $x_2\in Z_2\cap U_2$ respectively. Then $\theta_2\leq \theta_1$.
\end{lemma}
\proof
Since each chart $\theta_i$ is centered, we get an isomorphism $P_i \to \ol \M_{X,x}$, for any $x\in Z_i\cap U_i$, $i=1,2$. The composition
$\psi:P_1 \simeq \Gamma(U_1,\ol \M_X)\stackrel{\text{restr.}}{\longrightarrow}\Gamma(U_2,\ol \M_X)\simeq P_2$
is a surjection, inducing an isomorphism of $P_2$ with a localization of a face of $P_1$. 
Moreover, since for any $x\in U_2$ the  morphisms $\theta_{1,x}$ and $\theta_{2,x}\circ\psi$ induce the same map to $\ol\M_{X,x}$, there exists a homomorphim $\eta_x:P_1\to\O^{\times}_{U_2,x}$, such that $(\theta_{2,x}\circ\psi)\cdot\eta_x=\theta_{1,x}$.
Set $\varphi_x:=(\psi,\eta_x)$.
\endproof

\begin{lemma}\label{lem:directed charts on stratum}
	Let $\theta_U:P_U\to\Gamma(U,\M_X|_U)$ be a log chart, $V\subset U$ open and $p\in V$. 
	Assume that for any stratum $Z\subseteq X$ such that $V\cap Z\neq \emptyset$, we have that $V\cap \ol Z$ deformation retracts to $p$. 
	Then there exists a chart $\theta_V:P_V\to\Gamma(V,\M_X)$ on $V$ centered at $p$.
\end{lemma}
\proof
Let us consider the following diagram
\[	\begin{tikzpicture}[descr/.style={fill=white,inner 		
	sep=2.5pt}]
\matrix (m) [matrix of math nodes, row sep=3em,
column sep=3em]
{ P_U & \Gamma(U,\M_X)  & \Gamma(U,\ol\M_X)\\
	& \Gamma(V,\M_X)  & \Gamma(V,\ol\M_X)\\ };
\path[->,font=\scriptsize]
(m-1-1) edge node[above] {$\theta_U$} (m-1-2) 
(m-1-2) edge node[left] {$ $} (m-1-3)
(m-1-3) edge node[right] {$ $} (m-2-3)
(m-1-2) edge node[right] {$\text{restr.}$} (m-2-2) 
(m-2-2) edge node[left] {$ $} (m-2-3);
\end{tikzpicture},\]
Set $P_V:=\Gamma(V,\ol\M_X)$ and $\theta_U|_V:=\text{restr.}\circ\theta_U$.  Let
\[Q:=\{p\in P_U|\theta_U(p)|_V\in\Gamma(V,\O_X^{\times})\}.\]
Since for any stratum $Z\subseteq X$, $V\cap \ol Z$ deformation retracts to $p$, we have that $\Gamma(V,\ol \M_X)=\ol\M_{X,p}$. Then $P_V=P_U/Q$. Let us consider the exact sequence
\[0\to Q^{\gp}\to P_U^{\gp}\stackrel{\pi}{\to}P_U^{\gp}/Q^{\gp}\to 0.\]
Since $P_U^{\gp}/Q^{\gp}$ is torsion free, the sequence splits. Hence, there exists a map $\rho:P_U^{\gp}/Q^{\gp}\to P_U^{\gp}$ such that $\pi\circ\rho=\Id$. Then $\rho(P_U/Q)\subset P+Q^{\gp}$. Let
\begin{equation}
\begin{split}
\gamma:P_U+Q^{\gp}&\to \Gamma(V,\M_X)\\
(p,q)&\mapsto \theta_U(p)|_V\cdot (\theta_U(q)|_V)^{-1}\\
\end{split}.
\end{equation}
Set $\theta_V:=\gamma\circ\rho:P_V\to \Gamma(V,\M_X)$. By construction $\theta_V$ is centered at $p$.
\endproof	

\begin{proposition}\label{prop:existence directed log charts}
	There exists a finite set of directed log charts covering $(X,\M_X)$.
\end{proposition}
\proof
Let $(\theta'_i:P_i\to\Gamma(U_i,\M_X)_{i\in J})$ be a finite collection of log charts covering  $(X,\M_X)$. Let $(K, I_{\bullet})$ be a triangulation of $X$ \textit{adapted} to the stratification of $X$ induced by the log structure $\M_X$. This means that $K$ induces a triangulation $K_Z$ on each stratum $Z\subset X$. 
Up to refining the triangulation $K$ of $X$, via barycentring subdivisions, 
we can assume that for each vertex $v_i\in I_0$, there exists $i\in J$ such that $\text{Star}(v_i)\subseteq U_i$. By Lemma \ref{lem:directed charts on stratum},  for each vertex $v_i\in I_0$, we get a chart $\theta_{v_i}:P_{v_i}\to \Gamma(\text{Star}(v_i),\M_X)$ centered at $v_i$ (Definition \ref{def:centered charts}). 

Now, let $(v_{i_0},...,v_{i_k})\in I_{k}$. Without loss of generality, we assume $Z_{i_\mu}\subseteq \overline{Z}_{i_\nu}$, whenever $\mu<\nu$.
Consider
\[\bigcap^k_{j=0}\text{Star}(v_{i_j})=\text{Star}(\omega_{0\cdot\cdot\cdot k}),\]
where $\omega_{0\cdot\cdot\cdot k}$ is the minimal cell containing $v_{i_0},..,v_{i_k}$.
For each $l\in\{0,...,k\}$, let $(v_{i_0},..,\hat{v}_{i_l},...,v_{i_k})\in I_{k-1}$ obtained by  removing the element $v_{i_l}$. We have an inclusion map
\[\text{Star}(\omega_{0\cdot\cdot\cdot k})\hookrightarrow \text{Star}(\omega_{0\cdot\cdot\hat{l}\cdot\cdot k}).\]
By Lemma \ref{lem:equivalent charts}, we get \[\theta_{v_{i_{\nu}}}|_{\text{Star}(\omega_{0\cdot\cdot\cdot k})}\leq\theta_{v_{i_\mu}}|_{\text{Star}(\omega_{0\cdot\cdot\hat{l}\cdot\cdot k})},\] for any $\mu< \nu$, with $\nu\neq l$. 

Thus, for each $(v_{i_0},...,v_{i_k})\in I_k$, with $k\in\{0,1,2\}$, 
take 	\[\theta_{v_{i_k}}|_{\text{Star}(\omega_{0\cdot\cdot\cdot k})}:P_{v_{i_k}}\to \Gamma(\text{Star}(\omega_{0\cdot\cdot\cdot k}),\M_X).\] We remark that $P_{v_{i_k}}$ is the monoid with smallest rank among $(P_{v_{i_j}})_{j=0,...,k}$. Moreover,   for each $k\in\{1,2\}$ and $(v_{i_0},...,v_{i_k})\in I_k$, $l\in\{0,1,2\}$ and  $(w_{i_0},...,w_{i_{k-1}}):=(v_{i_0},..,\hat{v}_{i_l},...,v_{i_k})\in I_{k-1}$, take
\[\varphi^{w_{i_{k-1}}}_{v_{i_k}}:P_{w_{i_{k-1}}}\to P_{v_{i_k}}\oplus\O^{\times}_{\text{Star}(\omega_{0\cdot\cdot\cdot k})},\]
given by Lemma \ref{lem:equivalent charts}.

\endproof

	
	\bibliographystyle{acm}
	\bibliography{bibliography}

\begin{thebibliography}{10}

\bibitem{Raf}
{\sc Caputo, R.}
\newblock {\em Existence of a versal deformation for compact complex analytic
  spaces endowed with logarithmic structure}.
\newblock PhD thesis, Universität Hamburg, Von-Melle-Park 3, 20146 Hamburg,
  2019.
\newblock \url{http://ediss.sub.uni-hamburg.de/volltexte/2020/10273}.

\bibitem{DouVar}
{\sc Douady, A.}
\newblock Le probl\`eme des modules pour les vari\'{e}t\'{e}s analytiques
  complexes (d'apr\`es {M}asatake {K}uranishi).
\newblock In {\em S\'{e}minaire {B}ourbaki, {V}ol. 9}. Soc. Math. France,
  Paris, 1964, pp.~Exp. No. 277, 7--13.

\bibitem{DouThesis}
{\sc Douady, A.}
\newblock Le probl\`eme des modules pour les sous-espaces analytiques compacts
  d'un espace analytique donn\'{e}.
\newblock {\em Ann. Inst. Fourier (Grenoble) 16}, fasc. 1 (1966), 1--95.

\bibitem{DouVers}
{\sc Douady, A.}
\newblock Le probl\`eme des modules locaux pour les espaces {${\bf
  C}$}-analytiques compacts.
\newblock {\em Ann. Sci. \'{E}cole Norm. Sup. (4) 7\/} (1974), 569--602 (1975).

\bibitem{Felten}
{\sc Felten, S.}
\newblock Log smooth deformation theory via gerstenhaber algebras.
\newblock {\em manuscripta mathematica 167\/} (01 2022), 1--35.

\bibitem{Fischer1976}
{\sc Fischer, G.}
\newblock {\em Complex analytic geometry}.
\newblock Lecture Notes in Mathematics, Vol. 538. Springer-Verlag, Berlin-New
  York, 1976.

\bibitem{Flenner}
{\sc Flenner, H.}
\newblock \"{U}ber {D}eformationen homomorpher {A}bbildungen.
\newblock {\em Osnabr\"{u}cker {S}chriften zur {M}athematik\/} (1979).

\bibitem{Flenner1981}
{\sc Flenner, H.}
\newblock Ein kriterium für die offenheit der versalität.
\newblock {\em Mathematische Zeitschrift 178\/} (1981), 449--474.

\bibitem{ForKno}
{\sc Forster, O., and Knorr, K.}
\newblock {\em Konstruktion verseller {F}amilien kompakter komplexer
  {R}\"{a}ume}, vol.~705 of {\em Lecture Notes in Mathematics}.
\newblock Springer, Berlin, 1979.

\bibitem{GraLine}
{\sc Grauert, H.}
\newblock Analytische {F}aserungen \"{u}ber holomorph-vollst\"{a}ndigen
  {R}\"{a}umen.
\newblock {\em Math. Ann. 135\/} (1958), 263--273.

\bibitem{Gra}
{\sc Grauert, H.}
\newblock Der {S}atz von {K}uranishi f\"{u}r kompakte komplexe {R}\"{a}ume.
\newblock {\em Invent. Math. 25\/} (1974), 107--142.

\bibitem{GrossKansas}
{\sc Gross, M.}
\newblock {\em Tropical geometry and mirror symmetry}, vol.~114 of {\em CBMS
  Regional Conference Series in Mathematics}.
\newblock Published for the Conference Board of the Mathematical Sciences,
  Washington, DC; by the American Mathematical Society, Providence, RI, 2011.

\bibitem{GSGW}
{\sc Gross, M., and Siebert, B.}
\newblock Logarithmic {G}romov-{W}itten invariants.
\newblock {\em J. Amer. Math. Soc. 26}, 2 (2013), 451--510.

\bibitem{EGA43}
{\sc Grothendieck, A.}
\newblock \'{E}l\'{e}ments de g\'{e}om\'{e}trie alg\'{e}brique. {IV}. \'{E}tude
  locale des sch\'{e}mas et des morphismes de sch\'{e}mas. {III}.
\newblock {\em Inst. Hautes \'{E}tudes Sci. Publ. Math.}, 28 (1966), 255.

\bibitem{Gt}
{\sc Grothendieck, A.}
\newblock Techniques de construction et th\'{e}or\`emes d'existence en
  g\'{e}om\'{e}trie alg\'{e}brique. {IV}. {L}es sch\'{e}mas de {H}ilbert.
\newblock In {\em S\'{e}minaire {B}ourbaki, {V}ol. 6}. Soc. Math. France,
  Paris, 1995, pp.~Exp. No. 221, 249--276.

\bibitem{Hart}
{\sc Hartshorne, R.}
\newblock {\em Algebraic geometry}.
\newblock Springer-Verlag, New York-Heidelberg, 1977.
\newblock Graduate Texts in Mathematics, No. 52.

\bibitem{illusie2013}
{\sc Illusie, L., Nakayama, C., and Tsuji, T.}
\newblock On log flat descent.
\newblock {\em Proc. Japan Acad. Ser. A Math. Sci. 89}, 1 (01 2013), 1--5.

\bibitem{KatoF}
{\sc Kato, F.}
\newblock Log smooth deformation theory.
\newblock {\em Tohoku Math. J. (2) 48}, 3 (1996), 317--354.

\bibitem{KatoMa}
{\sc Kato, K.}
\newblock Logarithmic degeneration and {D}ieudonn\'{e} theory.
\newblock Available at
  \url{www.math.brown.edu/~abrmovic/LOGGEOM/Kato-Dieudonne.pdf}.

\bibitem{KatoFI}
{\sc Kato, K.}
\newblock Logarithmic structures of {F}ontaine-{I}llusie.
\newblock In {\em Algebraic analysis, geometry, and number theory ({B}altimore,
  {MD}, 1988)}. Johns Hopkins Univ. Press, Baltimore, MD, 1989, pp.~191--224.

\bibitem{KatNak}
{\sc Kato, K., and Nakayama, C.}
\newblock Log {B}etti cohomology, log \'{e}tale cohomology, and log de {R}ham
  cohomology of log schemes over {${\bf C}$}.
\newblock {\em Kodai Math. J. 22}, 2 (1999), 161--186.

\bibitem{kaup}
{\sc Kaup, L., Kaup, B., Barthel, G., and Bridgland, M.}
\newblock {\em Holomorphic Functions of Several Variables: An Introduction to
  the Fundamental Theory}.
\newblock De Gruyter studies in mathematics. W. de Gruyter, 1983.

\bibitem{KodNirSpr}
{\sc Kodaira, K., Nirenberg, L., and Spencer, D.~C.}
\newblock On the existence of deformations of complex analytic structures.
\newblock {\em Ann. of Math. (2) 68\/} (1958), 450--459.

\bibitem{Kur}
{\sc Kuranishi, M.}
\newblock On the locally complete families of complex analytic structures.
\newblock {\em Ann. of Math. (2) 75\/} (1962), 536--577.

\bibitem{Lepot}
{\sc Le~Potier, J.}
\newblock Le probl\`eme des modules locaux pour les espaces {${\bf
  C}$}-analytiques compacts (d'apr\`es {A}. {D}ouady et {J}. {H}ubbard).
\newblock 255--272. Lecture Notes in Math., Vol. 431.

\bibitem{Ogus18}
{\sc Ogus, A.}
\newblock {\em Lectures on Logarithmic Algebraic Geometry}.
\newblock Cambridge Studies in Advanced Mathematics. Cambridge University
  Press, 2018.

\bibitem{PalVers}
{\sc Palamodov, V.~P.}
\newblock Moduli in versal deformations of complex spaces.
\newblock In {\em Vari\'{e}t\'{e}s analytiques compactes ({C}olloq., {N}ice,
  1977)}, vol.~683 of {\em Lecture Notes in Math.} Springer, Berlin, 1978,
  pp.~74--115.

\bibitem{Pal}
{\sc Palamodov, V.~P.}
\newblock {\em Deformations of Complex Spaces, Several Complex Variables IV:
  Algebraic Aspects of Complex Analysis}, 1st~ed., vol.~10 of {\em
  Encyclopaedia of Mathematical Sciences}.
\newblock Springer-Verlag Berlin Heidelberg, 1990.

\bibitem{PourRel}
{\sc Pourcin, G.}
\newblock Th\'{e}or\`eme de {D}ouady au-dessus de {$S$}.
\newblock {\em Ann. Scuola Norm. Sup. Pisa (3) 23\/} (1969), 451--459.

\bibitem{Pourpriv}
{\sc Pourcin, G.}
\newblock Sous-espaces privil\'{e}gi\'{e}s d'un polycylindre.
\newblock {\em Ann. Inst. Fourier (Grenoble) 25}, 1 (1975), x--xi, 151--193.

\bibitem{ruddat2019local}
{\sc Ruddat, H.}
\newblock Local uniqueness of approximations and finite determinacy of log
  morphisms.
\newblock https://arxiv.org/abs/1812.02195, 2019.

\bibitem{RS}
{\sc Ruddat, H., and Siebert, B.}
\newblock Period integrals from wall structures via tropical cycles, canonical
  coordinates in mirror symmetry and analyticity of toric degenerations.
\newblock {\em Publications mathématiques de l'IHÉS 132}, 1 (2020), 1--82.

\bibitem{stacks-project}
{\sc {Stacks Project Authors}, T.}
\newblock \textit{Stacks Project}.
\newblock \url{https://stacks.math.columbia.edu/tag/00MP}.

\bibitem{Stie}
{\sc Stieber, H.}
\newblock {\em Existenz semiuniverseller {D}eformationen in der komplexen
  {A}nalysis}.
\newblock Aspects of Mathematics, D5. Friedr. Vieweg \& Sohn, Braunschweig,
  1988.

\end{thebibliography}

	

\end{document}